\documentclass[bj]{imsart}

\RequirePackage{amsthm,amsmath,amsfonts,amssymb}
\RequirePackage[numbers]{natbib}

\usepackage{amsmath, amssymb, tikz, epsfig, amsfonts, latexsym, mathrsfs, hyperref, graphicx, color, longtable, comment, bbm}

\startlocaldefs

\newcommand{\ccr}{C_c^+(\ol{\bbr}_0)}
\newcommand{\mbbo}{\mathbbm{1}}
\newcommand{\inv}{{-1}}
\newcommand{\bbr}{\mathbb{R}}
\newcommand{\probconv}{\stackrel{\prob}{\longrightarrow}}
\newcommand{\asconv}{\xrightarrow{{\rm a.s.}}}

\newcommand{\eqd}{\stackrel{{\rm d}}{=}}

\newcommand{\bbn}{\mathbb{N}}
\newcommand{\scrl}{\mathscr{L}}
\newcommand{\bbt}{\mathsf{T}}
\newcommand{\sfu}{\mathsf{u}}
\newcommand{\sfv}{\mathsf{v}}

\newcommand{\bdelta}{\boldsymbol{\delta}}
\newcommand{\scrm}{\mathscr{M}}
\newcommand{\dvague}{{\bf d}_{\rm vague}}
\newcommand{\uv}{{{\sf v}}}
\newcommand{\uu}{\sfu}
\newcommand{\Xu}{X(\uu)}

\newcommand{\wt}{\widetilde}
\newcommand{\wh}{\widehat}
\newcommand{\wtbfN}{\wt{\bf N}}
\newcommand{\ol}{\overline}

\newcommand{\bsymb}{\boldsymbol}
\newcommand{\poi}{\bsymb{\mathfrak P}}

\newcommand{\vep}{\varepsilon}

\newcommand{\h}{\hspace{.5cm}}

\newcommand{\mbbM}{\mathbb{M}}
\newcommand{\hlconv}{\xrightarrow{\mbbM_0}}
\newcommand{\mbbr}{\mathbb{R}}
\newcommand{\mbbn}{\mathbb{N}}
\newcommand{\uiid}{{\tiny{(\rm iid)}}}
\newcommand{\full}{{\tiny{(\rm full)}}}
\newcommand{\marked}{{\tiny{\rm(marked)}}}

\newcommand{\bbrz}{\bbr^\bbn_{\bf 0}}

\DeclareMathOperator{\prob}{\mathbb{P}}
\DeclareMathOperator{\dtv}{d}
\DeclareMathOperator{\exptn}{\mathbb{E}}

\DeclareMathOperator{\proj}{PROJ}

\newcommand{\beq}{\begin{equation}}
\newcommand{\eeq}{\end{equation}}
\newcommand{\alns}[1]{\begin{align*}#1\end{align*}}
\newcommand{\aln}[1]{\begin{align} #1 \end{align}}
\newcommand{\been}{\begin{enumerate}}
\newcommand{\een}{\end{enumerate}}
\newcommand{\norm}[1]{\| #1 \|}

\usepackage{cleveref}
\newtheorem{thm}{Theorem}[section]
\newtheorem{propn}[thm]{Proposition}

\newtheorem{lemma}[thm]{Lemma}

\theoremstyle{remark}
\newtheorem{remark}[thm]{Remark}

\theoremstyle{definition}
\newtheorem{defn}[thm]{Definition}
\newtheorem{fact}[thm]{{\bf Fact}}
\newtheorem{ass}[thm]{Assumptions}

\usepackage{todonotes}

\numberwithin{equation}{section}

\endlocaldefs

\allowdisplaybreaks

\begin{document}

\begin{frontmatter}
\title{Extreme positions of regularly varying branching random walk in random and time-inhomogeneous environment}
\runtitle{Extremes of BRW on BPRE}

\begin{aug}
\author[A]{\fnms{Ayan} \snm{Bhattacharya}\ead[label=e1, mark]{ayanbh@math.iitb.ac.in}}
\and
\author[B]{\fnms{Zbigniew} \snm{Palmowski}\ead[label=e2,mark]{zbigniew.palmowski@gmail.com}}
\address[A]{Department of Mathematics,
Indian Institute of Technology Bombay, Mumbai, India.
\printead{e1}}

\address[B]{Faculty of Pure and Applied Mathematics, Hugo Steinhaus Center\\
              Wroc\l aw University of Science and Technology\\
              Wyb. Wyspia\'nskiego 27, 50-370 Wroc\l aw, Poland,
\printead{e2}}
\end{aug}

\begin{abstract}
In this article, we consider a Branching Random Walk (BRW) on the real line. The genealogical structure is assumed to be given through a supercritical branching process in the i.i.d. environment and satisfies the Kesten-Stigum condition. The displacements coming from the same parent are assumed to have jointly regularly varying tails. Conditioned on the survival of the underlying genealogical tree, we prove that the appropriately normalized (normalization depends on the quenched size of the $n$-th generation) maximum among positions at the $n$-th generation converges weakly to a scale-mixture of Frech\'{e}t random variable. Furthermore, we derive the weak limit of the point processes composed of appropriately scaled positions at the $n$-th generation and show that the limit point process is a member of the randomly scaled scale-decorated Poisson point processes. Hence, an analog of the predictions by Brunet and Derrida \cite{brunet:derrida:2011} holds.  We have obtained an explicit description of the limit point process. This description captures the influence of the environment in the joint asymptotic behavior of the extreme positions. We show that the law of the clusters in the limit depends on the time-reversed environment (environment experienced by the rightmost/leftmost position). The asymptotic distribution of the normalized rightmost position is derived as a consequence. This approach (based on weak convergence of extremal processes) can not be adapted when the genealogical structure is given through a supercritical Branching Process in a time-Inhomogeneous Environment (BPIE) (due to lack of structural regularity in the genealogical structure). We provide a simple example where the point processes do not converge weakly. The tightness of the point processes holds in this example though the point processes do not (weakly) converge (due to having different subsequential weak limits). This phenomenon is not yet known in the literature.
\end{abstract}


\end{frontmatter}

\section{Introduction}

We consider a Branching Random Walk (BRW) on real line when law of the genealogical tree is given by a Branching Process  in the time-Inhomogeneous (random) Environment (BPIE) ${\bf Y} = ( Y_i : i \ge 0 )$ ($Y_i$ takes values in the set of probability measures on $\bbn_0$). To be more specific, the the process BRW starts with only one particle at the origin of the real line. We refer to it as the $0$-th generation of the process. Conditioned on the sequence ${\bf Y}$, each particle at the $(n-1)$-th generation dies, after producing an independent copy of the point process
\aln{
\scrl_n := \sum_{i =1}^{Z_1^{(n)}} \bdelta_{X_i}, \label{eq_defn_progeny_pp}
}
where law of $Z_1^{(n)}$ is given by $Y_{n-1}$ and $X_i$ denotes the {\it displacement} attached to the $i$-th child. Each newborn particle in the $n$-th generation are placed in the real line according to their {\it position} where position of a particle is given by its displacement translated by the position of its parent. If we run this procedure for long enough time, the resultant process of the positions might be called a BRW indexed by BPIE.

 In this article, we study the joint asymptotic behavior of the extreme positions at the $n$-th generation by establishing weak convergence of the point processes. Assuming the genealogical tree to be given by the supercritical BP in an i.i.d. environment, we show that the growth of the extreme positions depends on the quenched size of the genealogical tree up to the $n$-th generation.  We show that the limit point process is a randomly Scaled  scale-Decorated Poisson Point Process (SScDPP) which extends the predictions of Brunet and Derrida \cite{brunet:derrida:2011} to a time-inhomogeneous genealogical structure.  Due to the strong dependence structure among the positions, it is natural to expect that the extreme positions exhibit a clustering phenomenon.  We further show that the law of the clusters depends on the law of the environment reversed in time (see Theorems~\ref{thm_weak_pp_easy_version} and \ref{thm:main:thm}). To the best of our knowledge, this time reversal phenomenon is not reported yet in the literature on BRW. {\it The exchangeability property} of the environment (see proof of Proposition~\ref{propn:laplace:regular:point:process}) plays a crucial role in the proof of weak convergence of the point process. Although the principle of single-big-jump for regularly varying displacements makes the study of the extreme positions easier compared to the light-tailed displacements, the analog of our result is not yet known when the displacements are light-tailed.


 To underline the importance and significance of the i.i.d. environment in this paper, let us recall the main idea which stands behind the works on the extremes of the regularly varying BRW's. At the first step, one estimates the growth of the maxima of i.i.d. copies of the displacement where the size of the i.i.d. sample is the same as the size of the genealogical tree up to the $n$-the generation. The next step is technical and in this step, one divides the position at the $n$-th generation by the rate of growth and establishes weak convergence of the point processes composed of scaled-positions. The second step can be heuristically justified by the principle of single large displacement assuming the displacements to be i.i.d. It is tempting to check whether this picture is universal or not. To be more specific, if there is a supercritical branching process where this simple picture does not hold.
 We will provide a counterexample in Subsection~\ref{subsec_bpie} assuming the genealogical tree to be given by a BPIE.
 In particular, in Theorems~\ref{thm_tightness_rightmost_position} and \ref{thm_tighness_pp_simple_BPPE}, we will prove that
 above general idea cannot be applied for the case of he supercritical Branching Processes in the Periodic Environment (BPPE)
 with a constant period.
We shall show that the sequence of point processes (composed of scaled positions) does not converge weakly as the subsequential weak limits turn out to be different.  Still, this sequence of point processes
remains tight.
Because of periodicity, one may think that generation-dependent scaling might resolve this issue. But due to fundamental difference (dependence structure induced by the environment) in the distribution of the subsequential limits, the change in scaling is not enough to ensure the weak convergence of the point processes.
We conclude that the weak convergence of the point process may not be extended beyond the i.i.d. environment, which is assumed in this paper.

\subsection{Notations and facts from the branching process in i.i.d. environment}

The genealogical tree is denoted by ${\sf T} = ({\sf V}, {\sf E})$ where ${\sf V}$ and ${\sf E}$ denote the collection of all vertices and the edges. Throughout the paper, we shall follow the Harris-Ulam formalism to label the vertices of the genealogical tree. $\sfu $ and $\sfv$ denote typical vertices and $|\sfu|$ denotes the generation of the vertex $\sfu$.
The size of the $i$-th generation $\{ \sfu \in {\sf T} : |\sfu| = i\}$ will be denoted by $Z_i$ for $i \ge 0$.  The environment ${\bf Y}=(Y_i : i \ge 0)$ is a sequence of i.i.d. random variables
with the state space $\boldsymbol{\Theta}$ being a collection of probability measures on $\bbn_0$. We have assumed that $Z_0 = 1$.  A generic element of $\bsymb{\Theta}$ will be denoted by $\theta$ and $\exptn_\theta(\xi)$ will denote the expectation of the random variable $\xi$ distributed according to the law $\theta$.
We will use $\prob_{\bf Y}(\cdot)$ to denote the conditional probability $\prob(\cdot | {\bf Y}) $. Define
\begin{equation}\label{calS}
{\cal S} := \bigcap_{i \ge 1} \{ Z_i > 0\}.
\end{equation}
Let $q_e({\bf Y}) := \prob({\cal S} | {\bf Y})$ denote the probability of the survival conditioned on the environment
and $\prob^*(\cdot)=\prob(\cdot|{\cal S})$ denote the probability of an event given survival of the genealogical tree.
The conditional expectations associated to the conditional probabilities $\prob_{\bf Y}$ and $\prob^*$ will be denoted respectively by $\exptn_{\bf Y}$ and $\exptn^*$. $\mbbo_{{\sf A}}$ denotes the indicator function of the event ${\sf A}$. Define $\log^+(x) = \log (x \vee 1)$ and $\log^-(x) = \log (x^\inv \vee 1)$. We shall use $J_n \Rightarrow J$ to denote that $J_n$ converges weakly or in distribution to the random element $J$.

\begin{ass}
\label{ass:bpre:iidenv}
 $(Y_i : i \ge 0)$ is a sequence of i.i.d. random elements (taking values in the set of all probability measures on $\bbn_0$) which satisfy the following conditions:
\aln{
& \exptn \big| \log \prob (Z_1 > 1|Y_0) \big| < \infty,  \label{eq:ass:neg:mass:at:zero}\\
&  \exptn \Big( \log^{-} \exptn(Z_1|Y_0)  \Big) < \exptn \Big( \log^+ \exptn(Z_1 |Y_0) \Big) < \infty,  \label{eq:ass:supcritcal:iid:tree}  \\
& \mbox{ and } \exptn \Big( \frac{1}{\exptn( Z_1 |Y_0)} \exptn \big[ \mbbo_{\{Z_1 \ge 2\}} ~~ Z_1 \log Z_1 |Y_0\big]  \Big) < \infty. \label{eq:ass:ks:bpre}
}
\end{ass}

These assumptions imply that the underlying genealogical tree survives with the positive probability and grows indefinitely conditioned on the survival. Informally, \eqref{eq:ass:supcritcal:iid:tree} stands for {\it supercriticality} and \eqref{eq:ass:ks:bpre} is the analog of the Kesten-Stigum condition for BP in random environment. The implications of these assumptions are listed formally in the following lemma (see Theorem~1 in \cite{athreya:karlin:1971b}).

\begin{lemma}[Growth of supercritical BP in i.i.d. environment] \label{propn:supcrit:bp:iid:env}
Consider the sequence $(W_n := Z_n/ \pi_n : n \ge 0)$ where $Z_0 = \pi_0 =1$ and
\begin{equation}\label{firstpi}
\pi_n := \exptn_{{\bf Y}}(Z_n) = \prod_{i =0}^{n - 1} \exptn_{Y_i} (\xi )
\end{equation}
 for every $n \ge 1$. Then the random variables $(W_n : n \ge 1)$ form a non-negative martingale sequence with respect to the filtration $\big( \sigma\{Z_0,Z_1, Z_2, \ldots, Z_n;~~ {\bf Y}\} ~~: n \ge 0 \big)$
 and hence, there exists a random variable $W$ such that $\lim_{n \to \infty} W_n = W$, $\prob_{\bf Y}$-almost surely. If additionally the conditions \eqref{eq:ass:neg:mass:at:zero}, \eqref{eq:ass:supcritcal:iid:tree} and \eqref{eq:ass:ks:bpre} hold, then
\begin{enumerate}

\item $\exptn_{\bf Y}(W) = 1$ almost surely;

\item $1-q_e({\bf Y}) = \prob_{{\bf Y}}({\cal S}^c) = \prob_{{\bf Y}}(W = 0)$
almost surely, that is, the random variable $W$ is positive $\prob_{\bf Y}$-almost surely conditioned on the survival of the genealogical tree.
\end{enumerate}
\end{lemma}

\subsection{Asymptotic behavior of the extreme positions when displacements are i.i.d.}

Throughout the article, we shall assume that the displacements are identically distributed with regularly varying tails; that is,
\aln{
\prob(|X_1| > t) = t^{- \alpha}  {\rm L}(t) \mbox{ and } \lim_{t \to \infty} \frac{\prob(X_1 > t)}{\prob(|X_1| > t)} = p \in [0,1],  \label{eq:ass:marginal:regvar}
}
where ${\rm L}$ is a slowly varying function ($\lim_{t \to \infty} {\rm L}(tx)/ {\rm L}(t) = 1$) and $\alpha >0$. {\it For the sake of simplicity, we shall further assume in this subsection that the displacements are independently distributed and $\prob(Z_1 \ge 1) = 1$}. For a precise formulation, we need to develop more notations. Let ${\sf o}$ denote the root of  ${\sf T}$ and ${\sf o} \mapsto \uv$ denotes the geodesic path from the root ${\sf o}$ to the vertex $\uv$. The displacement attached to the vertex $\sfu \in {\sf T}$ will be denoted by $\Xu$ where $X({\sf o}) = 0$ almost surely. The position of the vertex $\uv \in {\sf T}$ is denoted by $S(\uv) = \sum_{\uu \in {\sf o} \mapsto \uv } \Xu$.  The collection $(S(\uv) : \uv \in {\sf T})$ is called {\it branching random walk}. In order to study the joint asymptotic behavior of the extreme positions (near-maxima), we will first investigate the asymptotic behavior of the maximum/rightmost position. Then,
in Subsection~\ref{subsec_pp_weak_iid_no_leaf}, we will use the growth of the maximum position to understand the joint behavior of the positions which has the same growth.

\subsubsection{Weak limit of the rightmost position}\label{subsec_rightmost_position}

The maximum/rightmost position in the $n$-th generation is denoted by $M_n = \max_{|\uv| = n} S(\uv)$. To understand the asymptotic behavior of $M_n$, one has to guess first a sequence $B_n$ (based on Lemma~\ref{propn:supcrit:bp:iid:env} and ``principle of single large displacement'' for the displacements) such that $(B_n^\inv M_n : n \ge 1)$ is a tight sequence of random variables. The next task would be to establish that $B_n^\inv M_n$ converges weakly to a random variable $M_*$ and identify the distribution of $M_*$. After the fundamental work \cite{durrett:1983}, it is natural to expect that $M_*$ is a scale-mixture of Fr\'{e}chet random variable.


\medskip

\noindent{\bf Guessing $B_n$.} As the genealogical tree does not have any leaves, its survival is certain.   Mathematically, the \textit{principle of single large displacement} means
\aln{
\prob \big( M_n > B_n x \big) \sim \prob \Big( \cup_{|\sfv| \le n} \{ X_\sfv > B_n x\} \Big) \mbox{ for every } x >0. \label{eq_principle_single_large_disp}
}
In words, a large displacement causes large values of the position of its descendants. One might get an annealed estimate of $\{\sfv \in \bbt: |\sfv| \le n\}$ (using i.i.d. structure of the environment), but that estimate can not be used for normalizing sequence to have a non-degenerate weak limit. So we turn to the quenched (conditioned on the environment) estimate $\pi_n$ of $\{|\sfv| \in \bbt : |\sfv| = n\}$, suggested in Lemma~\ref{propn:supcrit:bp:iid:env}. Under the assumptions~\ref{ass:bpre:iidenv}, it follows from SLLN that
\aln{
n^\inv \log \pi_n \asconv \exptn(\log \exptn(Z_1|Y_0)) > 0, \label{eq_geometric_growth_BPRE}
}
that is, $(\pi_n : n \ge 1)$ grows geometrically with high probability (w.h.p.). Thus, in our context, we may consider the ``quenched version'' of the principle of single large displacement (laws of the environment and displacements are assumed to be independent), that is,
\aln{
\prob \big( M_n > B_n x \big| {\bf Y} \big) \sim \prob \Big( \cup_{|\sfv| \le n} \{ X_\sfv > B_n x\} \big| {\bf Y} \Big) \label{eq_quenched_single_large_disp}
}
almost surely. From \eqref{eq_geometric_growth_BPRE}, it follows that
$
\exptn(\sum_{i = 1}^n Z_i | {\bf Y}) = \sum_{ i= 1}^n \pi_i \sim \pi_n
$
with high probability.
With appropriate (slightly more than second order) moment assumptions on the members of $\boldsymbol{\Theta}$, one can show that with high probability, the following asymptotics hold
$
\prob \Big( \cup_{|\sfv| \le n} \{ X_\sfv > B_n x \} | {\bf Y} \Big) \sim \pi_n \prob( X_1 > B_n x | {\bf Y}). 
$
 So, given the environment, this observation leads to the following definition of normalizing sequence
\aln{
B_n := \inf \{ s > 0 : \prob( |X_1| > s) < \pi_n^\inv \}.
\label{eq_defn_Bn}
}
If ${\rm L}$ is a constant function, then $B_n = \pi_n^{1/\alpha}$. If ${\rm L}$ is not a constant function,  $B_n = \pi_n^{1/\alpha} \wt{\rm L}(\pi_n)$ where $\wt{\rm L}$ is another slowly varying function (see Proposition~0.8(v) in \cite{resnick:1987}). We can show further that the weak limit of $B_n^\inv M_n$ exists and is a scale-mixture of Fr\'{e}chet random variable as presented in the next theorem. The following remark is important to understand the law of the scale-mixture.



\begin{remark}[Finiteness of the harmonic series of the quenched size of the supercritical BP in i.i.d. environment] \label{remark_harmonic_quenched_finite}
The harmonic series of the quenched sizes of generations in BPRE ${\bf Z}$ in environment ${\bf Y}$ is given by
$
1 + \sum_{n = 1}^\infty \pi_n^\inv = 1 + \sum_{n = 1}^\infty \Big[ \prod_{i = 0}^{n - 1} \exptn_{Y_i}(\xi) \Big]^\inv. 
$
Due to the exponential growth of $n \mapsto \pi_n$ (see \eqref{eq_geometric_growth_BPRE}), it is immediate that the harmonic series is almost surely finite.
\end{remark}

\begin{thm}[Weak limit of the rightmost position when displacements are i.i.d. and $\bbt$ does not have any leaf] \label{thm_maxpos_iid_disp_no_leaf}
Suppose that Assumptions~\ref{ass:bpre:iidenv} hold. Let ${\bf Y}'= (Y'_{i} : i \ge 0)$ be an independent copy of the environment sequence ${\bf Y}$ and hence, independent of $W$ in Lemma~\ref{propn:supcrit:bp:iid:env}.  Then there exists a random variable $M_*^{(iid)}$ such that $B_n^\inv M_n \Rightarrow M_*^{(iid)}$ such that
\aln{
&  \prob \Big( M_*^{(iid)} \le x \Big) = \exptn \Big[ \exp \Big\{ - W p x^{- \alpha} \Big( 1 +  \sum_{i =1}^\infty \Big[ \prod_{i = 0}^{n - 1} \exptn_{Y'_i}(\xi) \Big]^\inv \Big) \Big\} \Big]. \label{eq_maxima_iid_disp_no_leaf}
}
\end{thm}

\begin{remark}
 It is easy to check that if $\boldsymbol{\Theta}$ is singleton, then  $\prod_{i = 0}^{n - 1} \exptn_{Y'_j}(\xi) = m^n$ for every $n \ge 1$ where $m = \exptn(Z_1)$.   So Theorem~\ref{thm_maxpos_iid_disp_no_leaf} is an extension of the main result of \cite{durrett:1983} when the environment is time-inhomogeneous and can be derived easily from weak convergence of point processes (see Theorem~\ref{thm_weak_pp_easy_version}). The harmonic series in the exponent on the right-hand side of \eqref{eq_maxima_iid_disp_no_leaf} is not easy to guess from the results available in the literature. The most surprising element is that the random variable $W$ and the harmonic series are independent. Note that the largest displacement occurs in the last few generations (due to almost exponential growth of the genealogical structure) and each displacement has an equal probability to be the largest displacement. The random variable $W$ appears as the limit of the proportion of the eligible displacements to be the largest one and the harmonic series (which involves the number of their displacements) captures the contribution of the largest displacement to the maximum position.
\end{remark}


\subsubsection{Bulk behavior of the extreme positions} \label{subsec_pp_weak_iid_no_leaf}

The asymptotic behavior of the rightmost position provides only a partial picture of extreme positions as the rightmost position is usually not unique with positive probability. The extreme positions form clusters and may come with random multiplicities (the cluster size may have infinite mean even when the genealogical tree is a binary tree) due to the strong dependence structure among the positions. The most popular tool in extreme value theory is the point process which sheds some light on the joint asymptotic behavior of the extreme positions (comparable to the rightmost position). So we consider the point process composed of the scaled (scaling is same as that of the rightmost position) positions in the $n$-th generation and try to obtain the weak limit of these point processes. We start with basic facts and results on the point processes.

The asymptotic study of point processes is the most common tool in the extreme value theory to capture the bulk behavior of the extremes as it allows us to study the asymptotic behavior of the positions which are close to maxima.
Therefore we consider the following point process
\aln{
{\bf N}_n = \sum_{|\sfv| = n} \bdelta_{B_n^\inv S(\sfv)}
}
where $B_n$ is defined in \eqref{eq_defn_Bn}.

 Let $\scrm(\overline{\bbr}_0)$ be the space of all point measures on $\overline{\bbr}_0 : = [- \infty, \infty] \setminus \{0\}$.
We shall use $C_c^+(\overline{\bbr}_0)$ to denote the space of all non-negative bounded continuous functions on $\overline{\bbr}_0$ with compact support. It is known that vague convergence on $\scrm(\overline{\bbr}_0)$ is metrizable. The metric and topology induced by the vague convergence will be called {\it vague metric} and {\it vague topology} respectively. It can be shown that $\scrm(\overline{\bbr}_0)$ is a Polish space when equipped with the {\it vague topology} (see Proposition~3.17 in \cite{resnick:1987}). So, one can define the weak convergence of probability measures on $\scrm(\overline{\bbr}_0)$. To be mathematically more precise, a sequence of point processes $(\bsymb{\cal N}_n : n \ge 1)$  {\it converges weakly} to a point process $\bsymb{\cal N}_*$ in $\scrm(\overline{\bbr}_0)$ equipped with vague topology if and only if
\alns{
\lim_{n \to \infty} \exptn \Big[ \exp \Big\{ - \int f \dtv \bsymb{\cal N}_n \Big\} \Big] = \exptn \Big[ \exp \Big\{ - \int f \dtv \bsymb{\cal N} \Big\} \Big] \mbox{ for all } f \in C_c^+(\overline{\bbr}_0).
}
We refer to Proposition~3.19 in \cite{resnick:1987} for more details. In this paper, we shall denote it by $\bsymb{\cal N}_n \Rightarrow \bsymb{\cal N}_*$. $\exptn ( \exp \{ \int f \dtv \bsymb{\cal N}_*\})$ is called {\it Laplace functional} of $\bsymb{\cal N}_*$.  In this article, we shall establish that the Laplace functional of  ${\bf N}_n$ converges. Furthermore, we shall prove that the limit point process satisfies the heavy-tailed analog of prediction by Brunet and Derrida in \cite{brunet:derrida:2011}. To discuss it in a more formalized fashion, we need the following notations and definitions. We use $\mathscr{S}$ to denote the {\it scaling operator} on the space of all point measures such that ${\mathscr S}_a \big( \sum_{i \ge 1} \bdelta_{u_i} \big) = \sum_{i \ge 1} \bdelta_{a u_i} $ for every $a > 0$.

\begin{defn}[Randomly scaled scale-decorated Poisson point process]
A point process $\bsymb{\cal J}$ is called a {\it scale-decorated Poisson point process} (ScDPPP) with intensity measure $\kappa$ and scale decoration $\bsymb{\cal D}$  (denoted by $\bsymb{\cal J} \sim {\rm ScDPPP}(\kappa, \bsymb{\cal D})$) if there exist a Poisson random measure $\bsymb{\Lambda} = \sum_{i \ge 1} \bdelta_{\lambda_i}$ on $(0, \infty)$ with intensity measure $\kappa$ and a point process $\bsymb{\cal D} \in \scrm(\overline{\bbr}_0)$ such that
\aln{
\bsymb{\cal J} \eqd \sum_{i =1}^\infty \mathscr{S}_{\lambda_i} \bsymb{\cal D}_i,
}
where $(\bsymb{\cal D}_i : i \ge 1)$ is a collection of independent copies of $\bsymb{\cal D}$. Consider a positive random variable $U$ independent of the point process $\bsymb{\cal J}$, then the point process $\bsymb{\cal J}_* \eqd\mathscr{S}_U \bsymb{\cal J}$ is called a {\it randomly scaled scale decorated Poisson point process} (SScDPPP) and denoted by $\bsymb{\cal J}_* \sim {\rm SScDPPP}(\kappa, \bsymb{\cal D}, U)$.
\end{defn}

In \cite{bhattacharya:hazra:roy:2017}, it has been shown that every SScDPPP is a scale mixture of a strictly stable point process and the converse is also true.  It is well-known that the weak limit of the point processes in BRW is a member of SScDPPP when the genealogical structure is given by a supercritical Galton-Watson process. The most interesting feature of the limit point process (extremal process) is the law of the decoration point process. The following question comes naturally to mind: {\it How does the law of the genealogical structure influence the law of the clusters?} In most cases, it is difficult to answer the above question without a specific structural assumption on the dependence structure among the displacements coming from the same parent. In the case of displacements with exponentially decaying tails, the law of the decoration is yet unknown even when the genealogical tree is a binary tree. In Theorem~\ref{thm_weak_pp_easy_version} below, we shall provide the explicit description of the law of the clusters in terms of the law of the underlying BPRE.


Let $\nu_\alpha^{+}$ be the measure on $(0, \infty)$ such that $\nu_\alpha^{+}(x, \infty) = x^{- \alpha}$ for every $x > 0$.  Consider a Poisson point process $\sum_{l \ge 1} \bdelta_{\zeta_l}$ on $(0, \infty)$ with mean measure $\nu_\alpha^+$ which will be abbreviated as
\aln{
\sum_{l \ge 1} \bdelta_{\zeta_l} \sim {\rm PRM}(\nu_\alpha^{+}).
}
We consider a sequence of i.i.d. random variables $({\cal E}_l : l \ge 1)$ such that
\aln{
\prob({\cal E}_1 = 1) = 1 - \prob({\cal E}_1 = - 1) = p,
}
where $p$ is introduced in \eqref{eq:ass:marginal:regvar}. We assume that $({\cal E}_l : l \ge 1)$ is independent of the Poisson point process $\sum_{l \ge 1} \bdelta_{\zeta_l}$. Then, it is immediate to check that $\sum_{l \ge 1} \bdelta_{{\cal E}_l \zeta_l} \sim {\rm PRM} (\nu_\alpha)$ on $\overline{\bbr}_0$ where
\aln{
\nu_\alpha(x, \infty) = p \nu_\alpha^{+}(x, \infty) \mbox{ and }\nu_\alpha (- \infty, - x) = (1 - p) \nu_\alpha^{+}(x, \infty) \mbox{ for every } x > 0. \label{eq_defn_nu_alpha_measure}
}
 The later-mentioned Poisson point process is well-known in the literature on the classical extreme value theory (see \cite{davis:resnick:1985}).  We consider an independent copy ${\bf Y}' = (Y'_i : i \ge 0)$ of the environment sequence ${\bf Y}$ and define
\aln{
& {\rm C}_3({\bf Y}') =  1 + \sum_{i = 1}^\infty \Big[ \prod_{j = 0}^{i - 1} \exptn_{Y'_j}(\xi)  \Big]^\inv.
}
It follows from the Remark~\ref{remark_harmonic_quenched_finite} that ${\rm C}_3({\bf Y}') < \infty$ almost surely. Conditioned on ${\bf Y}'$, consider a collection $({\cal R}^{(l)} : l \ge 1)$ of independent copies of the positive integer-valued random variable ${\cal R}$. Here, the conditional probability mass function (p.m.f.)  of the random variable ${\bf Y}$ is given by
\aln{
& \prob \big( {\cal R} = r ~\big| ~ {\bf Y}' \big) = \frac{1}{{\rm C}_3({\bf Y}')} \Big[ \mbbo_{\{r = 1 \}} + \sum_{i = 1}^\infty \frac{\prob \big( Z_i = r ~| ~ {\bf Y}'_{i-1~:~0} \big)}{ \exptn \big( Z_i | {\bf Y}'_{i - 1~:~0} \big)} \Big]  \mbox{ for all } r \ge 1 \nonumber \\
& \hspace{1.5cm} \mbox{ where } {\bf Y}'_{i-1~:~0} = (Y'_{i - 1}, Y'_{i - 2}, \ldots, Y'_1, Y'_0) \mbox{ for every } i \ge 1.
}
Conditioned on ${\bf Y}'_{i- 1~: ~0}$, the random variable $Z_i$ denotes the size of the $i$-th generation of a BPRE such that the progeny distribution of a particle at the $j$-th generation is given by $Y'_{i - j}$ for all $j = 1, 2, \ldots, i$.
With these notations, we are ready to specify the distribution of the weak limit of ${\bf N}_n$.

\begin{thm}[Displacements are i.i.d. and genealogical tree does not have any leaf] \label{thm_weak_pp_easy_version}
There exists a point process $\widetilde{{\bf N}}_*^\uiid$ such that ${\bf N}_n \Rightarrow \widetilde{\bf N}_*^\uiid$ in the space $\scrm(\overline{\mbbr}_0)$ equipped with the vague topology. Furthermore,
\aln{
{\bf N}_*^{(iid)}  \eqd \sum_{l = 1}^\infty {\cal R}^{(l)} \bdelta_{[{\rm C}_3({\bf Y}') W]^{1/ \alpha} {\cal E}_l \zeta_l} &
\sim {\rm SScDPPP}(\nu_\alpha^+, {\cal R} \bdelta_{[{\rm C}_3({\bf Y}')]^{1/\alpha}{\cal E}},  W^{1/ \alpha}) \label{eq_rep_lim_pp_iid_disp}
}
and   for every $f \in \ccr$, $\exptn( \exp \{ - \int f \dtv {\bf N}_*\}) $ equals
\aln{
\exptn \Big[ \exp \Big\{ - W  \int_{\bbr_0} \nu_\alpha(\dtv x) \Big[ (1 - e^{- f(x)}) + \sum_{i = 1}^\infty \frac{  1 - \exptn ( e^{- f(x) Z_i} | {\bf Y}'_{i-1~:~0} )}{\exptn(Z_i | {\bf Y}'_{i - 1~:~0})} \Big] \Big\} \Big].   \label{eq_iid_disp_laplace_final}
}

\end{thm}

\begin{remark}
This theorem is a consequence of Theorem~\ref{thm:main:thm} in this paper. We refer to subsection~\ref{subsubsec_asymp_tail_indp} for a detailed discussion on the proof under a more general framework. The heuristics behind the point process convergence are as follows: Recall that the quenched size of the genealogical tree is of the order $\pi_n$ for large $n$. The principle of single large displacement and the definition of $B_n$ (see \eqref{eq_defn_Bn}) imply that the large (comparable to the magnitude of the maximal displacement) displacements occur in the last few ($n - o(n)$) generations with high probability due to exponential growth (see \eqref{eq_geometric_growth_BPRE}) of the underlying genealogical tree. The large displacements contribute to the position of their descendants in the $n$-th generation and make them comparable to the rightmost position. To establish weak convergence, one needs to keep counting the number of descendants (in the $n$-th generation) of the large displacements and their genealogical structure. The atoms of the limit ${\bf N}_*^{(iid)}$ correspond to the large displacements and clusters in the ${\bf N}_*^{(iid)}$ have the information on the genealogical structure among the extreme positions. The fraction of the number of large displacements is given by $W$ in the limit. The article aims to capture the influence of random environments on the law of the clusters. Note that extreme positions in the $n$-th generation see the environment in the reverse order. So, we have the time-reversed environment sequence in the limit. The exchangeability of the environment sequence turns out to be the key to the proof (see proof of Proposition~\ref{propn_laplace_limn}). Although time reversal of the environment in the limit is easy to anticipate, we are not aware of any other work on BRW formally establishing it.
\end{remark}

\begin{remark}[Genealogical tree is a supercritical Galton-Watson tree]
Assume that $\bsymb{\Theta}$ is singleton. Under the assumptions in \eqref{ass:bpre:iidenv}, $(Z_i : i \ge 0)$ turns out to be a supercritical Galton-Watson process that satisfies the Kesten-Stigum condition. So Theorem~\ref{thm_weak_pp_easy_version} reduces to the main result (Theorem~2.1) in \cite{bhattacharya:hazra:roy:2016} (we partially recover the main result of \cite{durrett:1983} as a consequence of Theorem~\ref{thm_maxpos_iid_disp_no_leaf}) when the progeny distribution does not put any mass at zero. As the necessary changes in the expression \eqref{eq_rep_lim_pp_iid_disp} and \eqref{eq_iid_disp_laplace_final} are straightforward, we omit the details.
\end{remark}

\begin{remark}
After generalizing the results on the asymptotic behavior of the extremes in the regularly varying BRW, it is tempting to believe that this approach (based on weak convergence of the point process) can be adapted to the supercritical BPIE. It is well-known in the literature that the asymptotic shape of the supercritical BPIE is not regular. In Subsection~\ref{subsec_bpie}, with the help of an simple example, we show that the sequence of point processes may not converge weakly for any BPIE. However, this example is not enough to conclude that weak convergence of point processes may not hold for a reversible ergodic Markovian environment. For a detailed discussion, we refer to the aforementioned subsection.
\end{remark}

\begin{remark}
Theorem~\ref{thm_maxpos_iid_disp_no_leaf} can be obtained as a consequence of Theorem~\ref{thm_weak_pp_easy_version} via the continuous mapping theorem. We derive Theorem~\ref{thm_weak_pp_easy_version} as a consequence of Theorem~\ref{thm:main:thm} though, that is derived under more general conditions. In a general framework (see Theorem~\ref{thm:main:thm}),  we assume that the displacements coming from the same parent have jointly regularly varying tail (see \eqref{eq_ass_disp_jtregvar} and, Subsection~\ref{subsubsec_jtregvar} and~\ref{subsec_regvar_polish_space} for a detailed discussion on extension of the concept of regular variation on punctured Polish space.) The concept of joint regular variation can accommodate independent displacements (or asymptotically tail-independent displacements) and fully dependent displacements as well. See Subsections~\ref{subsubsec_asymp_tail_indp}, \ref{subsubsec_lin_tail_dep} and \ref{subsubsec_full_tail_dep} for detailed discussion. The main issue there is that the SScDPPP representation for the weak limit is not explicit without any particular assumption on the dependence structure.
\end{remark}

\begin{remark}
Recently in \cite{xavier:2022}, the tightness of the appropriately shifted rightmost position has been obtained in the following framework. The displacements are assumed to be independently distributed according to the law of the standard Gaussian random variable. The genealogical tree is assumed to be given by an i.i.d. environment ${\bf Y}$ ($Y_i$'s are degenerated distributions and so can be interpreted as a random number) such that each particle at the $n$-th generation reproduces $Y_n$ particles in the next generation given the environment. In this case, the quenched size of the $n$-th generation becomes simple and is given by $\prod_{i = 1}^{n} Y_{i - 1}$. Suppose that the displacements are i.i.d. with regularly varying tails in the framework of the aforementioned reference. Then Theorem~\ref{thm_weak_pp_easy_version} can be used to conclude existence of a point process ${\bf N}_{**}$ such that ${\bf N}_n \Rightarrow {\bf N}_{**}$  where $B_n = (\prod_{i = 1}^n Y_{i - 1})^{1/ \alpha} \wt{\rm L}(\prod_{i = 1}^n Y_{i - 1})$. Due to the simplified form of quenched sizes, the Laplace functional of ${\bf N}_{**}$  equals
\alns{
\exptn \Big[ \exp \Big\{ - W \int_{\ol{\bbr}_0} \nu_\alpha(\dtv x) \Big[ (1 - e^{- f(x)}) + \sum_{i = 1}^\infty (\prod_{j = 1}^{i} Y'_{j - 1})^\inv \Big( 1 - \exp \Big\{ - \prod_{j = 1}^i { Y}'_{j - 1} f(x) \Big\} \Big) \Big] \Big\} \Big]
}
for every $f \in \ccr$. The SScDPPP representation of ${\bf N}_{**}$ turns out to be more interesting. Consider a sequence $({\cal R}_*^{(l)} : l \ge 1)$ of independent copies of the random variable ${\cal R}_*$ given ${\bf Y}'$ where
\alns{
\prob ({\cal R}_* = r) = [{\rm C}_4({\bf Y}')]^\inv \big( \mbbo_{\{ 0\}}(r) + \mbbo_{\bbn}(r) \big[ \prod_{j = 1}^r Y'_{j - 1} \big]^\inv \big) \mbox{ where } {\rm C}_4({\bf Y}') = 1 + \big[  \sum_{i = 1}^\infty \prod_{j = 1}^i Y'_j \big]^\inv.
}
Then it follows that ${\bf N}_{**} \eqd \sum_{l = 1}^\infty [\pi_{{\cal R}_*^{(l)}}({\bf Y}')]\bdelta_{[W {\rm C}_4({\bf Y}')]^{1/\alpha} {\cal E}_l \zeta_l}$. Furthermore, ${\bf N}_{**}$ admits a representation as ${\rm SScDPP}( \nu_\alpha^+, [\pi_{{\cal R}^*}({\bf Y}')] \bdelta_{[W {\rm C}_4({\bf Y}')]^{1/\alpha} {\cal E}}, W^{1/\alpha})$. Using the continuous mapping theorem, one may derive that $B_n^\inv M_n \Rightarrow M_{**}$ where   $\prob(M_{**} \le x) =  \exptn [ \exp\{ - W p x^{- \alpha} {\rm C}_4({\bf Y}') \}]$ for every $x \ge 0$.
\end{remark}




\subsection{Literature review}

BRW is considered to be an important model in statistical physics due to its connection to the models like Gaussian multiplicative chaos, Gaussian free field, first passage percolation, last passage percolation, scale-free percolation, and randomized algorithm. Therefore, a vast literature has emerged over the past few decades on the asymptotic behavior of the extreme positions. The pioneering works include \cite{hammersley1974}, \cite{kingman1975} and \cite{biggins1976}. We refer to \cite{shi:2015} for a recent review when the tails of the displacements decay exponentially.
Classically, BRW is indexed on a regular or Galton-Watson tree. We shall first discuss literature on the displacements admitting exponential moments. The order of fluctuation of the centered extreme positions has been obtained in \cite{bachmann:2000}, \cite{aidekon2013}, \cite{bramson:ding:zeitouni:2016}. Due to the prediction by Brunet and Derrida \cite{brunet:derrida:2011}, the bulk behavior of the extreme positions attracted many researchers. Recently, \cite{madaule:2011} established the weak convergence of the extremal process (seen from the minimal position) relying on the work of Maillard \cite{maillard:2013}. In \cite{subag:zeitouni:2014}, the invariance properties of the limit point processes have been studied. BRW in time-inhomogenous environment has been studied recently in \cite{fang:zeitouni:2012}, \cite{mallein:2015}, \cite{mallein:milos:2019}, \cite{ouimet:2018}, \cite{Liu1}, \cite{Liu2}, \cite{gao:13}, \cite{xavier:2022}. Most of these works consider the point process approach where the law of the displacements is also allowed to be time-dependent. Our work is closer to the framework of \cite{Liu1}, \cite{gao:13} and \cite{xavier:2022} where the progeny distribution is only allowed to be time-inhomogeneous.

The continuum analog of the BRW (when the tails of the displacements decay exponentially) is Branching Brownian Motion (BBM for later use). BBM is introduced for its maximal position that satisfies the famous FKPP equation. We refer to \cite{bovier:2017} and \cite{berestycki:2015} for a recent survey on BBM. Recently in \cite{hamel:nolen:2016}, the asymptotic properties of the KPP equation have been investigated assuming the underlying genealogical structure to be in a periodic environment. Relying on this study, the asymptotic behavior of the rightmost position in BBM has been studied in \cite{lubetzky:thornett:zeitouni:2018} when the genealogical structure is in a periodic environment (branching rate depends on the position of the particle). It has been shown that the rightmost position is tight when centered by its median and converges weakly to randomly shifted Gumbel distribution. The asymptotic order of the median is shown to be $O_{\prob} (\log t)$ (which is the same as that of a usual BBM) but the coefficient of $\log t$ becomes random. The asymptotic study of the point process is still open in this framework. The continuum analog of BRW is Branching Stable Motion (BSM for later use) when the displacements have a regularly varying tail with the index of regular variation in $(0, 2)$. This is still a blooming area of research which has been initiated in \cite{shiozawa:2008}, \cite{lalley:shao:2016} \cite{shiozawa:2021}. Most of the questions are open in this area.

We now discuss the works on the BRW when the displacements have so-called heavy-tails. The pioneering works include \cite{durrett:1979}, \cite{durrett:1983}, \cite{gantert2000}. Related recent works are \cite{maillard:2015}, \cite{berard2014}, \cite{gantert:hofelsauer:2018}, \cite{dyszewski:gantert:2020}, \cite{dyszewski:gantert:hofelsauer:2020}. The present article generalizes the previous weak convergence results obtained for the extremal process in the BRW with regularly varying displacements. In the literature of BRW, the main result Theorem~\ref{thm:main:thm} is the first one confirming the heavy-tailed analog of predictions by Brunet and Derrida when the genealogical structure is time-inhomogeneous. Although the proof of the main result is similar to the previous works, the explicit description of the limit is hard to guess given the literature. The time-reversal phenomenon is new though natural to anticipate. To the best of our knowledge, such a phenomenon is not explicitly obtained in any of the works which are available now. We believe that our analysis can be adopted to study the other BRWs with heavy-tailed displacements in a time-inhomogeneous environment.

\subsection*{Outline of the rest of the paper}

In Section~\ref{sec_main_results}, we shall recall necessary terminologies from $\mathbb{M}_0$ convergence of measures on the punctured Polish space, state the main result Theorem~\ref{thm:main:thm} along with its consequences. In Section~\ref{sec_proof_main_thm}, we discuss the main steps of Theorem~\ref{thm:main:thm} along with detailed proofs of Theorems~\ref{thm_tighness_pp_simple_BPPE} and \ref{thm_tightness_rightmost_position}. Proofs of the auxiliary results (in Section~\ref{sec_proof_main_thm}) can be found in the {\bf Appendix}.



\section{Main results} \label{sec_main_results}

We will first recall some basic notations and terminologies related to the $\mbbM_0$ convergence of measures on the punctured Polish space. These will be useful to mathematically specify the joint regular variation of the displacements coming from the same parent (see \eqref{eq_ass_disp_jtregvar}) using the regular variation of measures on the space of all real sequences. In Subsection~\ref{subsec_main_result}, after the necessary notations, we state the main result Theorem~\ref{thm:main:thm} concerning weak convergence of the extremal process. The consequences  are discussed in Subsection~\ref{subsec_conseq}. Finally, we shall discuss BRW indexed by BPIE in Subsection~\ref{subsec_bpie}.


\subsection{Regular variation on punctured Polish space and some results}

Let $\mathbb{S}$ be a Polish space equipped with metric $\rho$. The open ball of radius $r$ centered at $\bsymb{x}$ is denoted by ${\sf B}_r(\bsymb{x})$ and ${\cal B}(\mathbb{S})$ denotes the Borel $\sigma$-algebra generated by the open balls. Fix $\bsymb{s}_0 \in \mathbb{S}$ and let us consider the subspace ${\mathbb S}_0 = \mathbb{S} \setminus \{ \bsymb{s}_0\}$. ${\cal B}(\mathbb{S}_0)$ denotes the Borel $\sigma$-algebra on $\mathbb{S}_0$ induced by the subspace topology on $\mathbb{S}_0$. Suppose that ${\cal C}_b^+({\mathbb S}_0)$ denotes the class of all non-negative bounded continuous functions on $\mathbb{S}_0$ which vanishes in the open ball ${\sf B}_r(\bsymb{s}_0)$ for some $r > 0$. Let $\mathbb{M}({\mathbb S}_0)$ denote the space of all measures $\lambda$ on $\mathbb{S}_0$ such that $\int f \dtv \lambda < \infty$ for every $f \in {\cal C}_b^+(\mathbb{S}_0)$. A basic neighborhood of a measure $\lambda \in \mathbb{M}(\mathbb{S}_0)$ is given by $\{ \lambda' \in {\mathbb M}({\mathbb S}_0) : \max_{1 \le i \le k}|\int f_i \dtv \lambda' - \int f_i \dtv \lambda | < \epsilon\}$
 where $(f_i : 1 \le i \le k) \subset {\cal C}_b^+(\mathbb{S}_0)$.   So we can construct a sub-basis from the sets of the form
 \alns{
\{ \lambda \in \mathbb{M}(\mathbb{S}_0) : \int f \dtv \lambda \in {\sf G}\},~~ f \in {\cal C}_b^+ (\mathbb{S}_0) \mbox{ and } {\sf G} \mbox{ is an open subset of } [0, \infty).
 }
The topology on $\mathbb{M}(\mathbb{S}_0)$ induced by this sub-basis is called $\mathbb{M}_0$-topology. The convergence induced by this topology is called $\mathbb{M}_0$ convergence. If $(\lambda_n : n \ge 1)$ is a sequence of elements in $\mathbb{M}(\mathbb{S}_0)$ and converges to $\lambda$ in $\mathbb{M}_0$ topology, then we shall use $\lambda_n \hlconv \lambda$ to denote it. In this paper, we shall be mainly concerned with $\mathbb{S} = \mathbb{R}^\mbbn$ and $\bsymb{s}_0 = (0, 0 , 0 , \ldots)$. Of course there is a connection between the $\mbbM_0$-convergence and vague convergence, namely, $\mbbM_0$-convergence implies vague convergence that is, the compact sets in $\mbbM_0$- topology are compact in vague topology. We refer to Subsection~2.3 in \cite{lindskog:resnick:roy:2014} for a detailed and interesting discussion.

\subsubsection{Regular variation of measures on $\mbbr^\mbbn$ via $\mbbM_0$ topology} \label{subsec_regvar_polish_space}

Consider the metric $\rho_\infty$ on the space $\mbbr^\mbbn$ defined as $\rho_\infty({\bf x}, {\bf y}) = \sum_{i = 1}^\infty 2^{- i} (|x_i - y_i| \wedge 1) $ where ${\bf x} = (x_i : i \ge 1)$ and ${\bf y} = (y_i : i \ge 1)$ are two real sequences. It is clear that $(\mbbr^\mbbn, \rho_\infty)$ is a Polish space, and ${\cal B}(\mbbr^\mbbn)$ denotes the Borel $\sigma$-algebra generated by the open balls. Let ${\bf 0}_\infty = (0, 0, 0 , \ldots)$ be the origin of $\mbbr^\mbbn$. We can now use the $\mbbM_0$-convergence for the elements in $\mbbM(\mbbr^\mbbn \setminus \{ {\bf 0}_\infty\})$. Scalar multiplication plays a key role in the definition of regular variation. We can define scalar multiplication $\vartheta. {\bf x} = (\vartheta x_i : i \ge 1)$ for every sequence ${\bf x} = (x_i : i \ge 1)$ and $\vartheta > 0$ and it satisfies the following properties:
\begin{enumerate}
\item $(\vartheta, {\bf x}) \mapsto \vartheta . {\bf x} $ is continuous in $\vartheta$,
\item ${\bf 0}_\infty = \vartheta . {\bf 0}_\infty$ for every $\vartheta > 0$ and
\item $0 < \rho_\infty({\bf 0}_\infty, {\bf x}) \le  \rho_\infty({\bf 0}_\infty, \vartheta . {\bf x}) $ for all $\vartheta \ge 1$.
\end{enumerate}
We say a sequence $(t_n : n \ge 1)$ of positive scalars to be {\bf regularly varying} of index $\beta$ if $\lim_{n \to \infty} t_{\lfloor n \vartheta \rfloor}/t_n = \vartheta^{ \beta}$ for all $\vartheta > 0$ where $\lfloor x \rfloor$ denotes the largest integer less than $x$.  We denote it by $(t_n : n \ge 1) \in {\rm RV}_\beta$. If $\beta = 0$, we say $(t_n : n \ge 1)$ is {\bf slowly varying}.  We are now ready to define the regularly varying measure on $\bbrz = \bbr^\bbn \setminus \{{\bf 0}_\infty\}$.

\begin{defn} \label{defn_regvar_measure_on_bbrz}
A measure $\lambda \in \mbbM(\bbrz)$ is said to be regularly varying if there exists a non-null measure $\varsigma \in \mbbM(\bbrz)$ and a sequence $(t_n : n \ge 1) \in {\rm RV}_\beta$  for some $\beta > 0$ such that $t_n \lambda(n \cdot) \hlconv \varsigma(\cdot)$ where $\lambda(n . {\sf A}) = \lambda (\{n. {\bf x} : {\bf x} \in {\sf A}\})$ for every ${\sf A} \in {\cal B}(\bbrz)$.
\end{defn}

Note that there are many equivalent definitions of the regular variation and we refer to Subsection~3.2 in \cite{lindskog:resnick:roy:2014} for a detailed discussion on that. According to Theorem~3.1 in the aforementioned reference, the limit measure $\varsigma$ in definition~\ref{defn_regvar_measure_on_bbrz} must satisfy
\aln{
\varsigma(\vartheta. {\sf A}) = \vartheta^{-\beta} \varsigma({\sf A}) \mbox{ for every } \vartheta > 0 \mbox{ and for every } {\sf A} \in {\cal B}(\bbrz). \label{eq_homogeneity_prop_lim_measure}
}
That is, the measure $\varsigma$ satisfies a homogeneity property. To stress this, we shall use $\lambda \in {\rm RV}_{- \beta}(\bbrz, ~\varsigma)$ for a regularly varying measure $\lambda$ on $\bbrz$.

%

\subsubsection{Assumptions on displacements} \label{subsubsec_jtregvar}

We can use the notations and the terminologies developed in the previous subsections to specify the dependence structure of the displacements coming from the same parent.
\begin{ass}
\begin{enumerate}
\item Displacements are independent of the branching mechanism. The displacements ${\bf X} = (X_i : i \ge 1)$ are independent of $(Z_1^{(n)} : n \ge 1)$ in the definition of $\scrl_n$ (see \eqref{eq_defn_progeny_pp}).

\item The displacements are identically distributed with regularly varying tails, that is \eqref{eq:ass:marginal:regvar} holds for every $X_i$.

\item ${\bf X}$ is jointly regularly varying that is,
\aln{
\prob ({\bf X} \in \cdot) \in {\rm RV}_{- \alpha} (\bbrz, \bsymb{\nu}). \label{eq_ass_disp_jtregvar}
}

\end{enumerate}
\end{ass}
Under these assumptions on the displacements, we shall scale the positions at the $n$-th generation by $B_n$ (defined in \eqref{eq_defn_Bn}) which is guessed based on the i.i.d. displacements. Due to the joint regular variation, we can not have $\nu_\alpha$ in the limit but some measure related to $\bsymb{\nu}$.  If the BPRE satisfies the assumptions~\ref{ass:bpre:iidenv}, then it can be shown that
\aln{
& \pi_n \prob \Big( B_n^\inv {\bf X} \in \cdot \Big| {\bf Y}_{0~:~n - 1} \Big) \hlconv \bsymb{\nu}_*(\cdot) \mbox{ almost surely }   \label{eq_mzero_conv_crucial}\\
& \mbox{ where } \bsymb{\nu}_* = \bsymb{\nu}( \cdot ) /\bsymb{\nu} ({\sf E}_1) \mbox{ and } {\sf E}_1 = \big( \big\{ (x_i : i \ge 1) : |x_1| > 1 \big\} \big). \nonumber
}
We shall briefly discuss how the $\mbbM_0$-convergence in \eqref{eq_mzero_conv_crucial} follows from the assumptions on the displacements. We first observe that $\pi_n \prob(B_n^\inv {\bf X} \in {\sf E}_1) = \pi_n \prob(|X_1| > B_n) \to 1$ for almost all environment ${\bf Y}$ due to geometric growth of $\pi_n$ (see \eqref{eq_geometric_growth_BPRE}) and definition of $B_n$ (see \eqref{eq_defn_Bn}). It is easy to check that ${\sf E}_1$ is bounded away from ${\bf 0}_\infty$ that is, ${\bf 0}_\infty \notin {\rm cl}({\sf E}_1)$ where ${\rm cl}({\sf E}_1)$ denotes the closure of the set ${\sf E}_1$. It can be shown (using the absolute continuity of $\bsymb{\nu}$ induced by the homogeneity property of $\bsymb{\nu}$) that ${\bsymb \nu}(\partial {\sf E}_1)  = 0$ where $\partial {\sf E}_1$ denotes the boundary of the set ${\sf E_1}$.  This observation along with Theorem~3.1 in \cite{lindskog:resnick:roy:2014} (equivalence of definitions (ii) and (iii) of regularly varying measure on $\bbrz$) leads us to the following observation
\aln{
\lim_{n \to \infty} \prob \Big( B_n^\inv {\bf X} \in {\sf C} \Big| {\bf Y} \Big) / \prob \Big( B_n^\inv {\bf X} \in {\sf E}_1 \Big| {\bf Y} \Big) = \bsymb{\nu}({\sf C})/ \bsymb{\nu}({\sf E}_1) \mbox{ for every } {\sf C} \in {\cal B}(\bbrz)
}
for almost all environments. Hence, the almost sure $\mbbM_0$-convergence in \eqref{eq_mzero_conv_crucial} follows.

\subsection{Weak convergence of the rightmost position and extremal process} \label{subsec_main_result}

We shall first introduce some random variables and random processes (constructed on the same probability space $(\Omega, {\cal F}, \prob)$) which will be used to write down the weak limit of ${\bf N}_n$.

\begin{enumerate}
\let\myenumi\theenumi
\let\mylabelenumi\labelenumi
\renewcommand{\theenumi}{N\myenumi}
\renewcommand{\labelenumi}{{\rm (\theenumi)}}


\item ${\bf Y}' = (Y'_i : i \ge 0)$ is an independent copy of environment ${\bf Y}$.

\item $[i~:~j]$ denotes the set $\{i, i + 1, i +2, \ldots, j\}$ for $j > i  \ge 0$ and the set
$\{i, i -1, i -2, \ldots, j\}$ for $i > j  \ge 0$.  We define ${\bf Y}'_{i~:~0} = (Y'_i, Y'_{i -1}, Y'_{i -2}, \ldots, Y'_0)$ for $i \ge 1$ and ${\bf Y}'_{0~:~0} = Y'_0$.

\item $Z_1^{(+)}$ denotes the random variable $Z_1$ conditioned to be positive that is, $ \prob(Z_1^{(+)} = k) = \prob(Z_1 = k | Z_1 \ge 1)$.

\item The law of $\wt{Z}_i$ conditioned on the environment segment ${\bf Y}'_{n-1~:~0}=(Y'_{n-1}, Y'_{n - 2}\ldots, Y'_0)$ is  the law of the number of descendants in the $n$-th generation of a particle in the $(n - i)$-th generation for every $i \in [1:n-1]$.

\item $\wt{ Z}_i^{(+)}$ denotes the random variable $\wt{ Z}_i$ conditioned to be positive that is, $\prob(\wt{Z}_i^{(+)} = k | {\bf Y}'_{n - 1~:~0}) = \prob \big( \wt{ Z}_i = k | (Y'_{i -1}, Y'_{i -2}, \ldots, Y'_1, Y'_0); ~~ \{ \wt{ Z}_i \ge 1\} \big)$  for all  $ k \in \bbn$.

\item Conditioned on the environment ${\bf Y}'_{n - 1~:~0}$, $(\wt{Z}_i^{(+,k)} : k \ge 1)$ denotes a collection of independent copies of the random variable $\wt{ Z}_i^{(+)}$  for every $i \in [1~:~n-1]$.

\item ${\rm Pow}({\sf B}) := \{{\sf C} : {\sf C} \subset {\sf B}\}$ is the power set of the set ${\sf B}$.

\item  Consider a Poisson random measure $\poi$ on the space $\bbr^\bbn_{{\bf 0}}$ with the mean measure $\bsymb{\nu}_*$ (defined in \eqref{eq_mzero_conv_crucial}) such that
\aln{
\poi := \sum_{l \ge 1} \bdelta_{\bsymb{\zeta}^{(l)}} = \sum_{l =1}^\infty \bdelta_{(\zeta^{(l)}_1, \zeta_{2}^{(l)}, \zeta_{3}^{(l)}, \ldots )}. \label{eq_poisson_rm}
}
We assume that $\poi$ is independent of the random variable $W$ and environment sequence ${\bf Y}'$.

\item  Let $\{{\bf 0}_i\}$ be the origin of $\bbr^i$ and $\bbn_{\#}^i = \{0, 1, 2, \ldots, \infty\}^i \setminus \{{\bf 0}_i\} $ for every $i \ge 1$. Conditioned on ${\bf Y}'$, consider a $\bigcup_{i =1}^\infty \big\{ \{i \}\times \bbn_{\#}^i \big\}$-valued random variable  $(V, {\bf R}) = (V, (R_l : 1 \le l \le V))$ with the following conditional probability mass function
\aln{
& \prob_{{\bf Y}'}(V= v; ~ {\bf R} = {\bf r}) = \frac{1}{{\rm C}_1({\bf Y}')} \Big[ \sum_{i = 1}^\infty \frac{\prob(Z_1 = v | Y'_i)}{ \exptn(Z_{i + 1} | {\bf Y}'_{i~:~0})} \prod_{m = 1}^v \prob(Z_i = r_m | {\bf Y}'_{i-1~:~0}) \nonumber \\
& \hspace{2cm} + [\exptn(Z_1 | Y'_0)]^\inv \prob(Z_1 = v | Y'_0) \prod_{m = 1}^v \mbbo_{\{ r_m = 1\}} \Big]  \mbox{ for all } {\bf r} \in  \bbn_{\#}^v \mbox{ and } v \ge 1 \label{eq_cond_pmf_cluster} \\
& \mbox{ where } {\rm C}_1({\bf Y}') = [\exptn(Z_1 | Y'_0)]^\inv + \sum_{i = 1}^\infty \big[ \exptn( Z_{i + 1} | {\bf Y}'_{i~:~0}) \big]^\inv \Big[ \prob(Z_1 \ge 1 | Y'_i)  \nonumber \\
& \hspace{2cm} - \sum_{v = 1}^\infty \prob(Z_1 = v | Y'_i ) [ \prob(Z_i = 0| {\bf Y}'_{i - 1~:~0} )]^v \Big].
}
It follows immediately from Remark~\ref{remark_harmonic_quenched_finite} that ${\rm C}_1({\bf Y}') < \infty$ almost surely.  The random vector $(V, {\bf R})$ is constructed independently of the random variable $W$ and the Poisson process $\poi$.

\item Given the environment ${\bf Y}'$, consider a collection  $((V_l, {\bf R}^{(l)}) = (V_l , (R^{(l)}_k : 1 \le k \le V_l)): l \ge 1)$ of independent copies of $(V, {\bf R})$ which are also independent of $W$ and $\poi$.

\end{enumerate}

With these random variables and notations, we are now ready to present the main result of this article.
\begin{thm} \label{thm:main:thm}
If the assumptions \eqref{eq:ass:supcritcal:iid:tree}, \eqref{eq:ass:neg:mass:at:zero}, \eqref{eq:ass:ks:bpre}, \eqref{eq:ass:marginal:regvar} and \eqref{eq_ass_disp_jtregvar} hold, then there exists an ${\rm SScDPPP}$ ${\bf N}_*$ such that conditioned on the survival ${\cal S}$,  ${\bf N}_n \Rightarrow  {\bf N}_*$ in the space $\scrm(\ol{\bbr}_0)$ equipped with vague topology. Moreover, we have
\aln{
 {\bf N}_* \eqd  \sum_{l= 1}^\infty \sum_{k =1}^{V_l} R^{(l)}_k \bdelta_{( {\rm C}_1({\bf Y}') W)^{1/\alpha} \zeta^{(l)}_{k}}
\label{characterization}}
  and  for
$f\in \ccr $,  the Laplace functional $\exptn^* [\exp \{ - {\bf N}_*(f) \}]$ can be written as follows
\aln{
 & \exptn^* \Bigg[ \exp \Bigg\{ - W \int_{\bbr^\bbn_{{\bf 0}} }
\bsymb{\nu}_*
(\dtv {\bf x})  \Bigg[  \sum_{i =1}^{\infty} [ \exptn(Z_{i +1} | {\bf Y}'_{i~:~0}) ]^\inv \prob(Z_1 \ge 1  | Y'_i)  \exptn \Big[ \sum_{{\sf B} \in {\rm Pow}([1:~ Z_1^{(+)}]) \setminus \{\emptyset\}}   \nonumber \\
& \hspace{1cm}\big[ \prob \big( Z_i \ge 1 \big| {\bf Y}'_{i - 1~:~0} \big) \big]^{|{\sf B}|} \big[ \prob \big( Z_i = 0 | {\bf Y}'_{i - 1~:~0} \big) \big]^{Z_1^{(+)} - |{\sf B}|}  \Big( 1 - \exp \Big\{ - \sum_{k \in {\sf B}} \wt{ Z}_i^{(+,k)} f(x_k) \Big\} \Big) \Big| {\bf Y}'_{i ~:~0}  \Big] \nonumber \\
& \hspace{2cm} +  [\exptn(Z_1 | Y'_0)]^\inv \prob(Z_1 \ge 1 | Y'_0) \exptn \Big( 1 - \exp \Big\{ - \sum_{k = 1}^{Z_1^{(+)}} f(x_k) \Big\} \Big| Y'_0\Big)   \Bigg] \Bigg\} \Bigg]. \label{eq:final:expression:laplace:functional:main:thm}
}
\end{thm}

\begin{remark}[Clusters of extremes in the limit and quenched law of the branching process.]
It is clear from \eqref{characterization} that ${\bf N}_*$ is a Cox cluster process. The Poisson atoms $(\bsymb{\zeta}_l : l \ge 1)$ appear because the collection $({\bf X}_{\sfu} : \sfu \in {\sf T})$ are i.i.d. where ${\bf X}_{\sfu}$ denotes the displacement vector attached to the children of $\sfu$. One large displacement vector (principle of ``a bunch of large displacements'') causes multiple large positions and a strong dependence structure among the positions. This dependence structure induces multiplicities ($((R^{(l)}_k : 1 \le k \le V_l) : l \ge 1)$) of the Poisson atoms $(\bsymb{\zeta}^{(l)} : l \ge 1)$ where the law of the multiplicities is given by the quenched law of the generation-sizes in an independent environment. Due to the i.i.d. structure, the law of the environment sequence does not change even when reversed in time (this fact has been used in the proof of Proposition~\ref{propn:laplace:regular:point:process}). The most interesting point here is that the law of the clusters involves the time-reversed environment. The main reason is that we are looking at the large positions at the $n$-th generation and then looking genealogically backward to their most recent common ancestor with the large displacements. We would like now to explain the appearance of the random variable $W$ and ${\bf Y}'$, and their independence in the limit ${\bf N}_*$. It follows from Lemma ~\ref{propn:supcrit:bp:iid:env} and SLLN that the sizes of the generations grow geometrically with a high probability conditioned on the survival of the tree. If we combine this fact with the principle of ``a bunch of large displacements'', it follows that the large displacements can occur in the last few ($o(n)$ with high probability) generations. So we cut the tree and ignore all the displacements except the last few generations. Due to cutting the genealogical tree very close to the generation $n$, we create a forest containing i.i.d. trees with displacements and the number of trees in the forest is comparable to the size of the $n$-th generation of the original genealogical tree. The random variable $W$ appears from the number of trees in the forest after appropriate normalization. Thanks again to the i.i.d. structure of the environment sequence, the environment corresponding to the first $n - o(n)$ generations of the genealogical tree is independent of the environment of the forest. Hence the independence of $W$ and ${\bf Y}'$.
\end{remark}

\begin{remark}[${\bf N}_*$ as an SScDPPP]
We first note that it is not possible to obtain an explicit series representation of ${\bf N}_*$ as an SScDPPP due lack of information on $\bsymb{\nu}_*$. To verify that ${\bf N}_*$ is an SScDPPP, one has to borrow tools like scaled-Laplace functional which has been introduced in Section~3 of \cite{bhattacharya:hazra:roy:2017}. In this paper, we shall not do it as one can use the homogeneity property of the measure $\bsymb{\nu}_*$ (see \eqref{eq_homogeneity_prop_lim_measure}) to establish the claim following the proof of Theorem~2.6 in the aforementioned reference.
\end{remark}

We now turn to the asymptotic behavior of maximum position $M_n = \max_{|\uv| = n} S(\uv)$. The existence of the weak limit $M_*$ of $B_n^\inv M_n$ is an immediate consequence of the weak convergence of point processes due to
\aln{
\prob^* (M_* \le x) = \lim_{n \to \infty} \prob^* (B_n^\inv M_n \le x )= \lim_{n \to \infty} \prob^* \big( {\bf N}_n (x, \infty) = 0 \big) = \prob^*({\bf N}_* (x, \infty) = 0).
}
As ${\bf N}_*$ is an SScDPPP, it follows immediately from equation (4.3) in \cite{bhattacharya:2018b} (see  (Prop2) in Proposition~3.2 in \cite{bhattacharya:hazra:roy:2017}) that $M_*$ is scale-mixture of Frech\'{e}t random variable. That is, $M_* = Q. \Phi_\alpha$ where $\Phi_\alpha$ denotes the Frech\'{e}t-$\alpha$ random variable and $Q$ is a positive random variable independent of $\Phi_\alpha$. It is well known that the random variable $Q$ preserves the features of the underlying genealogical structure and the dependence structure among the displacements. Although $Q$ is the most important random variable, it is known to be a challenging problem (especially when the displacements have an exponentially decaying tail) to disintegrate the law of $Q$. As an explicit representation for the limit ${\bf N}_*$ is obtained in Theorem~\ref{thm:main:thm}, a detailed description of the law of $Q$  can be extracted and is given in the next theorem. As a preparation of the next result,  let us define ${\cal G}_0 = [- \infty, 1]$, ${\cal G}_1 = (1, \infty]$ and
\aln{
{\cal H}^{(t)}_{i_1, i_2, \ldots, i_t} : = \prod_{j = 1}^t {\cal G}_{i_j} \times \prod_{j'= t + 1}^\infty \bbr \mbox{ where } i_j \in \{0, 1\} \mbox{ for all } j = 1, 2, \ldots, t. \label{eq_defn_calH_maxpos}
}

\begin{thm}[Weak limit of the appropriately normalized rightmost position] \label{thm_maxpos_general}
Suppose that the assumptions of Theorem~\ref{thm:main:thm} hold. Then conditioned on the event ${\cal S}$, $B_n^\inv M_n \Rightarrow M_*$ such that for every $x > 0$, $\prob(M_* \le x) = \exptn( e^{- x^{- \alpha} Q})$ and
\alns{
Q \eqd & W \sum_{j = 0}^\infty [\exptn(Z_{j + 1} | {\bf Y}'_{j : 0})]^\inv \sum_{v = 1}^\infty \prob(Z_1 = v | Y'_j) \\
& \h \h \sum_{k = 1}^v \big[ 1 - (\prob[Z_j = 0 | {\bf Y}'_{j - 1 : 0}])^k  \big] \sum_{\substack{i_1, i_2, \ldots, i_v  \in \{0,1\} \\ i_1 + i_2 + \ldots + i_v = k}}  \bsymb{\nu}_*({\cal H}^{(v)}_{i_1, i_2, \ldots, i_v}). \label{eq_max_dep_final_form}
}
\end{thm}
As the theorem can be proved following the lines of the proof of Corollary~2.7 in \cite{bhattacharya:hazra:roy:2017} with appropriate modifications, we will not include a detailed proof in this paper. We only focus on underlying the main differences.

\begin{remark}[The genealogical tree does not have a leaf]
Suppose $\prob(Z_1 \ge 1) = 1$ implying survival of the genealogical tree be certain. The weak limit $\wh{\bf N}_*$ and its Laplace functional admit simplified forms. Given the environment sequence ${\bf Y}'$, consider a collection $((\wh{\bf R}^{(l)}, \wh{V}_l) : l \ge 1)$ of independent copies of the random vector $(\wh{\bf R}, \wh{V})$ with conditional p.m.f.
\aln{
& \prob_{{\bf Y}'} (\wh{V} = v; ~ \wh{\bf R} = {\bf r}) = \frac{{\rm C}_1({\bf Y}')}{{\rm C}_2({\bf Y}')} \prob_{{\bf Y}'} \big( V= v; {\bf R} = {\bf r} \big) \mbox{ for all } {\bf r} \in \bbn^v \mbox{ and } v \ge 1 \nonumber \\
&  \mbox{ where } {\rm C}_2({\bf Y}') = 1 +  \sum_{i = 1}^\infty [\exptn(Z_{i + 1} | {\bf Y}'_{i ~: ~0})]^\inv.
}
For every $f \in \ccr$, the Laplace functional $\exptn[\exp \{ - \wh{\bf N}_*(f)\}]$ equals
\aln{
& \exptn \Big[ \exp \Big\{ - W \int_{\bbr_{{\bf 0}}^{\bbn}} \bsymb{\nu}_*(\dtv {\bf x}) \Big[ [\exptn(Z_1 | Y'_0)]^\inv \exptn \big( 1 - \exp \big\{ - \sum_{k = 1}^{Z_1} f(x_k) \big\} | Y'_0 \big)  \nonumber \\
& \hspace{2cm} +  \sum_{i = 1}^\infty [\exptn(Z_{i + 1} | {\bf Y}'_{i~:~0} )]^\inv \exptn \big[ 1 - \exp \big\{ - \sum_{k = 1}^{Z_1} \wt{Z}_i^{(k)} f(x_k) \big\}\big|~ {\bf Y}'_{i~:~0} \big] \Big] \Big\} \Big]  \label{eq_mult_reg_var_no_leaf_laplace}
}
and $\wh{\bf N}_* \eqd \sum_{l = 1}^\infty \sum_{k = 1}^{\wh{V}_l} \wh{R}_k^{(l)} \bdelta_{[W {\rm C}_2({\bf Y}')]^{1/\alpha} \zeta^{(l)}_k}$. As a consequence, $B_n^\inv M_n \Rightarrow \wh{M}_*$ where $\prob(\wh{M}_* \le x)$ equals
\alns{
\exptn \Big[ \exp \Big\{ - W x^{-\alpha} \sum_{i = 0}^\infty [\exptn(Z_{i + 1} | {\bf Y}'_{i~:~0})]^\inv \sum_{v = 1}^\infty \prob(Z_1 = v | Y'_i) \sum_{\substack{i_1, i_2, \ldots, i_v \in \{0, 1\} \\ i_1 + i_2 + \ldots + i_v \ge 1}} \bsymb{\nu}_*({\cal H}^{(v)}_{i_1, i_2, \ldots, i_v}) \Big\} \Big]
}
where ${\cal H}_{i_1, i_2, \ldots, i_v}^{(v)}$ is introduced in \eqref{eq_defn_calH_maxpos}.
\end{remark}


\subsection{Consequences} \label{subsec_conseq}

Here, we mention some of the consequences of the Theorem~\ref{thm:main:thm} where an explicit description of the measure ${\bsymb \nu}$ is available. We first address the case where the displacements are asymptotically tail-independent. The weak limit of the point processes is an extension of Theorem~\ref{thm_maxpos_iid_disp_no_leaf} allowing the genealogical tree to have leaves. Then we address the case where the displacements are (asymptotically) fully tail-dependent and derive the weak limit of ${\bf N}_n$ in Theorem~\ref{thm_fully_tail_dependent} along with the SScDPPP representation. Finally, we consider the displacements from a simple linear process and obtain an explicit description of the weak limit of ${\bf N}_n$.

\subsubsection{Asymptotically tail-independent displacements} \label{subsubsec_asymp_tail_indp}

We shall call the displacements ${\bf X} = (X_i : i \ge 1)$ to be {\it asymptotically tail-independent} if there exists a sequence $(t_n : n \ge 1)$ such that
\aln{
t_n \prob \big( n^\inv {\bf X} \in \cdot \big) \hlconv \sum_{i = 1}^\infty {\bsymb \tau}_i \mbox{ where } {\bsymb \tau}_i = \bigotimes_{j = 1}^{i - 1} \bdelta_{0} \otimes \nu_\alpha \bigotimes_{j = i + 1}^\infty \bdelta_0 \mbox{ for all } i \ge 1. \label{eq_defn_tail_independent_disp}
}
That is, ${\bsymb \tau}_i$ concentrates on the $i$-th axis. See page~195 in  \cite{resnick:2007} (for finite dimensional asymptotically tail-independence) and references therein; and Subsection~4.5.1 in \cite{lindskog:resnick:roy:2014} for the infinite-dimensional analogue. Independent displacements are asymptotically tail-independent. It turns out that some dependent random variables may also turn out to be asymptotically tail-independent. As the joint regular variation (\eqref{eq_ass_disp_jtregvar}) can accommodate the asymptotically tail-independent random variables with ${\bsymb \nu}$ given by the right hand side of \eqref{eq_defn_tail_independent_disp}, Theorem~\ref{thm:main:thm} can be used to conclude existence of weak limit $\wt{\bf N}_*^{(iid)}$ of ${\bf N}_n$. Due to an explicit description of ${\bsymb \nu}$, one can obtain a simplified expression of the Laplace functional given by
\aln{
\exptn^* \Big[ \exp \Big\{ - W \int_{\bbr_0} \nu_\alpha(\dtv x) \Big[ (1 - e^{- f(x)}) +  \sum_{i = 1}^\infty \frac{\prob \big( Z_i \ge 1 | {\bf Y}'_{i -1~:~0} \big)}{ \exptn(Z_i | {\bf Y}'_{i - 1~:~ 0})} \exptn \Big( 1 - e^{f(x) Z_i^{(+)}} \Big| {\bf Y}'_{i - 1~:~0} \Big) \Big] \Big\} \Big].   \label{eq_final_laplace_iid_disp}
}
This can be thought of as an extension of Theorem~\ref{thm_maxpos_iid_disp_no_leaf} allowing displacements to be dependent and satisfying \eqref{eq_defn_tail_independent_disp}, and genealogical tree to have leaves.
It is possible to write down a series representation of the point process $\wt{\bf N}_*^{(iid)}$ but it involves some new notations.  Recall the random processes ${\bsymb {\cal P}}, {\bf Y}'$ and $({\cal E}_l : i \ge 1)$ from Subsection~\ref{subsec_main_result} which are independent of the martingale limit $W$ (introduced in Lemma~\ref{propn:supcrit:bp:iid:env}). Conditioned on the sequence ${\bf Y}'$, consider a collection of i.i.d. random variables $({\cal R}^{(l)} : l \ge 1)$ such that
\aln{
& \prob \big( {\cal R} = r | {\bf Y}' \big) = \frac{1}{{\rm C}_3({\bf Y}')} \Big[ \mbbo_{\{r = 1\}} +  \sum_{i = 1}^\infty \frac{\prob(Z_i = r | {\bf Y}'_{i-1~:~0})}{ \exptn(Z_i | {\bf Y}'_{i - 1~:~0})} \Big] \mbox{ for every } r \ge 1 \\
&   \mbox{ where } {\rm C}_3({\bf Y}') = 1 +  \sum_{i = 1}^\infty \frac{\prob(Z_i \ge 1 | {\bf Y}'_{i - 1~:~0})}{\exptn(Z_i | {\bf Y}'_{i - 1~:~0})} \mbox{ is the normalizing constant}.
}
Then, it can be shown that
\aln{
\wt{\bf N}_*^{(iid)} \eqd \sum_{l = 1}^\infty {\cal R}^{(l)} \bdelta_{[W{\rm C}_3({\bf Y}')]^{1/\alpha} {\cal E}_l \zeta_l} & = {\mathscr S}_{[W{\rm C}_3({\bf Y}')]^{1/\alpha} } \sum_{l = 1}^\infty {\mathscr S}_{\zeta_l} \big( {\cal R}^{(l)} \bdelta_{{\cal E}_l} \big) \nonumber \\
& \sim {\rm SScDPPP} \big(~ \nu_\alpha^+, {\cal R} \bdelta_{[{\rm C}_3({\bf Y}')]^{1/\alpha} {\cal E}}, ~ W^{1/\alpha}~ \big). \label{eq_final_series_rep_iid_disp}
}
This representation turns out to be immensely helpful to obtain an explicit expression form for the distribution of the weak limit $\wt{M}_*^{(iid)}$ (weak convergence of the rightmost position follows from the weak convergence of the point processes) as
\aln{
\lim_{n \to \infty} \prob^* \big( M_n \le x \big) = \prob^* \big( \wt{M}_*^{(iid)} \le x \big) = \exptn^* \big[ e^{- W p x^{- \alpha} {\rm C}_3({\bf Y}')}  \big].  \label{eq_iid_disp_max_pos_distn_with_leaves}
}
Note that Theorem~\ref{thm_maxpos_iid_disp_no_leaf} and \ref{thm_weak_pp_easy_version} can be derived from \eqref{eq_iid_disp_max_pos_distn_with_leaves} and \eqref{eq_final_laplace_iid_disp}, \eqref{eq_final_series_rep_iid_disp} respectively.

\begin{remark}[Joint asymptotic distribution of the rightmost and leftmost position]
As a consequence of the tail-balancing conditions, the asymptotic joint distribution of the maximum and the minimum position can be obtained from point process convergence via continuous mapping theorem. Let $\ol{{\cal M}}_n = \min_{|\uv| = n} S(\uv)$ and $\ol{\cal M}^{(iid)}_*$ be the weak limit of $\ol{\cal M}_n$. Then it follows that for every $x, y > 0$, $\lim_{n \to \infty} \prob^* \Big( \{ M_n \le B_n x \} \cap \{ \ol{ \cal M}_n > - B_n y\} \Big) $ equals
\alns{
 \prob^* \Big( \{\wt{ M}_*^{(iid)} \le x\} \cap \{ \ol{\cal M}_*^{(iid)} > -  y\} \Big)  = \exptn^* \Big[\exp \big\{  - W {\rm C}_3({\bf Y}') (p x^{- \alpha} - (1 - p) y^{- \alpha}) \big\}\Big].
}
using the fact that conditioned on $(W, {\bf Y}')$, the atoms of the point process ${\bf N}_*^{(iid)}$ are coming from a Poisson random measure.
\end{remark}

\begin{remark}[Joint asymptotic distribution of the first two of large positions $M_n , M_n^{(2)}$ and first gap] \label{remark_joint_distn_two_max}
From the weak convergence of $\wt{\bf N}_n$ it follows immediately that  $B_n^\inv(M_n , M_n^{(2)}) \Rightarrow (M_*^{\tiny{(\rm iid)}}, M_{*,2}^{\tiny{(\rm iid)}})$ under $\prob^*$. So it is a natural question if we can write down the joint distribution of $ (M_*^{\tiny{(\rm iid)}}, M_{*,2}^{\tiny{(\rm iid)}})$ explicitly. Consider $0 < y < x $. Then it follows immediately that $\prob^*(\{M_*^{\tiny{{\rm (iid)}}} \le x\} \cap \{ M_{*,2}^{\tiny{\rm (iid)}} \le y\} ) = \prob^* (M_*^{\rm (iid)} \le y)  + \prob^* (\{ y < M_*^{\rm (iid)} \le x \} \cap \{ M_{*, 2}^{\rm (iid)} \le y \})$.  It is easy to check that $ \prob^* ( \{ y < M_*^{\rm (iid)} \le x\} \cap \{ M_{*, 2}^{\rm (iid)} \le y\}) = \exptn^* [\prob^*( \{ {\bf N_*}(x, \infty) = 0\} \cap \{ {\bf N}_*(y, x] = 1\} | W, {\bf Y}')]$. It is easy to check that $\{ {\bf N}_* (x, \infty) = 0 \} =\{ \sum_{l \ge 1} \bdelta_{{\cal E}_l \zeta_l} ([W {\rm C}_3({\bf Y}')]^{-1/\alpha} x, \infty)\} = 0$ and $\{ {\bf N}_*(y, x] = 1\} = \{ \sum_{l \ge 1} \bdelta_{{\cal R}^{(l)}, [W {\rm C}_3({\bf Y}')]^{- 1/\alpha} {\cal E}_l \zeta_l} (\{1\} \times ( y ,x] = 1 \}$. Conditioned on $W$ and ${\bf Y}'$, $\sum_{l \ge 1} \bdelta_{{\cal E}_l \zeta_l}$ and $\sum_{l \ge 1} \bdelta_{{\cal R}^{(l)}, {\cal E}_l \zeta_l}$ are Poisson point process and so, we can use the independence along with some algebraic adjustments to conclude that
\aln{
& \exptn^* [\prob^*( \{ {\bf N_*}(x, \infty) = 0\} \cap \{ {\bf N}_*(y, x] = 1\} | W, {\bf Y}')] \nonumber \\
& =p \exptn^* \Big[\Big( 1 +  \sum_{i = 1}^\infty \frac{\prob(Z_i = 1 | {\bf Y}'_{i - 1~:~0})}{\exptn(Z_i | {\bf Y}'_{i-1~:~0})} \Big) \exp \Big\{ - W  p \Big[{\rm C}_3({\bf Y}') x^{- \alpha}  \nonumber \\
& \hspace{2cm} +  W(y^{- \alpha} - x^{- \alpha}) \Big( 1 + \sum_{i = 0}^\infty \frac{\prob(Z_i = 1 | {\bf Y}'_{i- 1~:~0})}{ \exptn(Z_i | {\bf Y}'_{i-1~:~0})} \Big] \Big) \Big\} W(y^{- \alpha} - x^{- \alpha})  \Big]. \label{eq_joint_distn_asymp_first_two_maxima}
}
An important function of $B_n^\inv (M_n^{(1)}, M_n^{(2)})$ is the gap statistic $B_n^\inv (M_n^{(1)} - M_n^{(2)})$.  From the formula \eqref{eq_joint_distn_asymp_first_two_maxima}, it is formally possible to derive the asymptotic distribution of the gap statistic. It is difficult to obtain an explicit description of the law of the gap statistic. In the last decade, gap statistics has gained its importance in the literature on statistical physics (see \cite{majumdar:pal:2014} for a review and the references therein).  In the case of a weak dependence structure among the extremes, the normalized gap between the first two large positions may turn out to be absolutely continuous (with respect to the Lebesgue measure).
If the atoms in the limiting point process have multiplicities (as in our case), the probability of the gap statistic equal to $0$ must be positive. So the law of the gap statistic is not absolutely continuous.   Therefore, the probability of the gap being zero turns out to be interesting to investigate in its own right. It will be interesting to derive the conditional density of the gap statistic given the gap statistic is positive as the gap statistic might put mass at zero. There is not much literature on the gap statistic which may help to answer these questions.
\end{remark}

\subsubsection{Asymptotically fully tail-dependent displacements} \label{subsubsec_full_tail_dep}

Suppose that the displacements coming from the same parent are equal that is, ${\scrl}^{(n)} \eqd Z_1^{(n)} \bdelta_{X_1}$ where $\prob(X_1 \in \cdot ) \in {\rm RV}_{-\alpha}(\bbr_0, \nu_\alpha)$; see page~196 in \cite{resnick:2007} for applications and motivations. Consider the sequence ${\bf 1}_\infty$ whose all components equal $1$. Then it follows immediately that
\aln{
\prob( {\bf X} \in \cdot) \in {\rm RV}_\alpha ( \bbr^\bbn_{{\bf 0}}, \bsymb{\nu}_\full) \mbox{ where } {\bsymb \nu}_\full(\cdot) = \nu_\alpha(\{ x \in \bbr: x {\bf 1}_\infty \in \cdot\}).  \label{eq_hl_limit_full_dep}
}
Displacements satisfying \eqref{eq_hl_limit_full_dep} are called {\it asymptotically fully tail-dependent}. Recall the random variable $W$, random processes ${\bsymb {\cal P}}$, ${\bf Y}'$, $(({\bf R}^{(l)}, V_l) : l \ge 1)$ and $({\cal E}_l : l \ge 1)$. The following theorem can easily be derived from Theorem~\ref{thm:main:thm} replacing $\bsymb{\nu}_*$ by $\bsymb{\nu}_\full$ after some algebraic adjustments.

\begin{thm}\label{thm_fully_tail_dependent}
Suppose that the assumptions~\ref{ass:bpre:iidenv}, \ref{eq:ass:marginal:regvar} and \eqref{eq_hl_limit_full_dep} hold. Then there exists a point process ${\bf N}_*^\full$ such that conditioned on ${\cal S}$, ${\bf N}_n \Rightarrow {\bf N}_*^\full$ in the space $(\scrm(\overline{\bbr}_0), \dvague)$. Furthermore, we have
\aln{
{\bf N}_*^\full &\eqd \sum_{l =1}^\infty \Big( \sum_{k =1}^{V_l} R_k^{(l)} \Big) \bdelta_{[{\rm C}_1({\bf Y}') W]^{1/ \alpha} \zeta_l {\cal E}_l}  = \mathscr{S}_{[{\rm C}_1({\bf Y}') W]^{1/\alpha} } \sum_{l= 1}^\infty \mathscr{S}_{\zeta_l} \Big(\bdelta_{{\cal E}_l}  \sum_{k =1}^{V_l} R_k^{(l)}  \Big) \nonumber \\
& \sim {\rm SScDPPP}(\nu_\alpha^+, \bdelta_{[{\rm C}_1({\bf Y}') ]^{1/\alpha}{\cal E}} \sum_{k =1}^V R_k, W^{1/\alpha})
}
and ${\bf N}_*^\full$ has the following Laplace functional
\aln{
& \exptn^* \Big[ \exp \Big\{ - W \int_{\bbr_0} \nu_\alpha(\dtv x) \Big[ \sum_{i =1}^\infty \frac{\prob(Z_1 \ge 1 | Y'_i)}{\exptn(Z_{i + 1} | {\bf Y}'_{i~:~0})}  \exptn  \Big( \sum_{{\sf B} \in {\rm Pow}([1~:~ Z_1^{(+)}])} \big[  \prob( Z_i \ge 1 | {\bf Y}'_{i-1~:~0} \big)\big]^{|{\sf B}|} \nonumber \\
& \h \h \big[ \prob(Z_i = 0 | {\bf Y}'_{i-1~:~0}) \big]^{Z_1^{(+)} - |{\sf B}|} \big[ 1 - e^{- f(x) \sum_{k \in {\sf B}} Z_i^{(+, k)}}\big]  ~\Big|~ {\bf Y}'_{i~:~0}  \Big) \nonumber \\
& \h \h \h  + \frac{\prob(Z_1 \ge 1 | Y'_0)}{ \exptn(Z_1 | Y'_0)} \exptn \big[ 1 - e^{- Z_1^{(+)} f(x)} \big| Y'_0 \big] \Big] \Big\} \Big] \label{eq_lim_Laplace_full_dependence}
}
for every $f \in \ccr$. As a consequence, $M_n \Rightarrow M_*^\full$ under $\prob^*$ and for every $x > 0$, we have following expression for $\prob^*(M_*^\full \le x)$
\aln{
 & \exptn^* \Big[ \exp \Big\{ - W p x^{- \alpha} \Big[  \sum_{i = 1}^\infty \frac{\prob(Z_1 \ge 1 | { Y}'_i)}{ \exptn(Z_{i + 1} |  {\bf Y}'_{i~:~0})} \exptn \Big( 1  \nonumber \\
 & \h \h - \big[ \prob(Z_i = 0 | {\bf Y}'_{i-1~:~0}) \big]^{Z_1^{(+)}}  \Big| {\bf Y}'_{i~:~0} \Big) + \frac{\prob(Z_1 \ge 1 | Y'_0)}{\exptn(Z_1 | Y'_0)}\Big] \Big\} \Big].
}
\end{thm}

\begin{remark}[Joint asymptotic distribution of the first two large positions]
Fix $0 < y \le x$. We can use the above result along with the observations in Remark~\ref{remark_joint_distn_two_max} to obtain
\aln{
& \lim_{n \to \infty} \prob^* \Big( \{ B_n^\inv M_n^{(1)} \le x\} \cap \{ B_n^\inv M_n^{(2)} \le y \}  \Big) = \prob^* \Big(  \{ M_*^{\tiny{({\rm full})}} \le x\} \cap \{ M_{*, 2}^{\tiny{({\rm full})}} \le y\} \Big) \nonumber \\
& = \prob^*(M_*^{\rm (full)} \le y) + \prob^* (\{ y < M_*^{\rm (full)} \le x\} \cap \{ M_{*, 2}^{\rm (full)} \le y\}) \nonumber \\
& = \prob^*( M_*^{\rm (full)} \le y) +   p (y^{- \alpha} - x^{- \alpha})\exptn^* \Big[ W {\rm C}_1({\bf Y}')\prob^* \Big(  M_*^{\tiny{({\rm full})}} \le x | {\bf Y}', W \Big) \prob(V= 1; R = 1 | {\bf Y}')  \nonumber\\
& \h \h   \exp \Big\{ - p W {\rm C}_1({\bf Y}') (y^{- \alpha} - x^{- \alpha}) \prob(R = 1; V= 1 | {\bf Y}') \Big\} \Big].
}
\end{remark}

\subsubsection{Linearly dependent displacements and multivariate regular variation} \label{subsubsec_lin_tail_dep}

Assume that there exists an integer ${\rm K} > 1$ such that $\prob(Z_1 \le {\rm K}) = 1$. Then, we can use the multivariate regular variation to model the tail dependence among the displacements coming from the same parent. We refer to Subsection~6.1.4 in \cite{resnick:2007} for a detailed discussion on the multivariate regular variation. We shall discuss a particular case where the displacements are linearly dependent.

Let $L_0$, $L_1$ and $L_2$ be i.i.d. random variables such that $\prob(L_0 \in \cdot) \in {\rm RV}_{- \alpha}((0, \infty), \nu_\alpha^+)$. Define the random variables
\aln{
X_1 = L_1 + \phi L_0 \mbox{ and } X_2 = L_2 + \phi L_1.
}
Consider a random variable $Y$ such that $p_1 = \prob(Y = 1) = 1 - \prob( Y = 2)$. Let $(Y_i : i \ge 0)$ be a collection of independent copies of the random variable $Y$. Conditioned on the sequence ${\bf Y} = (Y_i : i \ge 0)$, each particle in the $n$-th generation reproduces an independent copy of the point process $\sum_{i = 1}^{Y_n} \bdelta_{X_i}$.  This yields a BRW in an i.i.d. environment where displacements coming from the same parent are linearly dependent.

Note that $(L_0, L_1, L_2) \in {\rm RV}_{- \alpha}([0, \infty)^3 \setminus \{{\bf 0}_3\}, \bsymb{\tau})$ where ${\bsymb \tau}(\cdot) = (\nu_\alpha^+ \otimes \bdelta_0 \otimes \bdelta_0 + \bdelta_0 \otimes \nu_\alpha^+ \otimes \bdelta_0 + \bdelta_0 \otimes \bdelta_0 \otimes \nu_\alpha^+)(\cdot)$. Consider the map $\wt{\rm T}_\phi : [0, \infty)^3 \to [0, \infty)^2$ such that $\wt{\rm T}_\phi(l_1, l_2, l_3) = (l_2 + \phi l_1, l_3 + \phi l_2)$.
One can use continuous mapping theorem (Theorem~2.3 in \cite{lindskog:resnick:roy:2014}) for the $\mbbM_0$ convergence to conclude that
\aln{
& \prob( (X_1, X_2) \in \cdot) \in {\rm RV}_{- \alpha} ([0, \infty)^2 \setminus \{ {\bf 0}_2\}, \bsymb{\tau} \circ \wt{\rm T}_\phi^\inv) \mbox{ where for all } {\sf A} \in {\cal B}([0, \infty)^2 \setminus \{ {\bf 0}_2\}),\nonumber \\
&  \bsymb{\tau} \circ \wt{\rm T}_\phi^\inv ({\sf A}) = \bsymb{\tau}(\{ (l_1, l_2, l_3) \in [0, \infty)^3 \setminus \{{\bf 0}_3\}: (l_2 + \phi l_1, l_3 + \phi l_2) \in {\sf A}\}).
}
We can easily obtain the weak limit ${\bf N}_*^{\rm (lin)}$ of the point process in this case as a consequence of Theorem~\ref{thm:main:thm} replacing $\bsymb{\nu}_*$ by $\bsymb{\tau} \circ \wt{\rm T}_\phi^\inv$. It is straightforward to check from \eqref{eq_mult_reg_var_no_leaf_laplace} that the Laplace functional of the weak limit ${\bf N}_*^{{\tiny {\rm (lin)}}}$ equals
\aln{
& \exptn \Big[ \exp \Big\{ - W \int_{\bbr^2_{{\bf 0}}} \bsymb{\tau} \circ \wt{\rm T}_\phi^\inv (\dtv {\bf x}) \Big[ (Y'_0)^\inv \Big( 1 - \exp \Big\{ - \sum_{k = 1}^{Y'_0} f(x_k) \Big\} \Big)  \nonumber \\
& \hspace{1cm} +  \sum_{i = 1}^\infty (\prod_{j = 0}^{i} Y'_j)^\inv  \Big( 1 - \exp \big\{ - (\prod_{j = 0}^{i - 1} Y'_j) \sum_{k = 1}^{Y'_i} f(x_k) \big\}  \Big)  \Big] \Big\} \Big].
}
Let $\sum_{l \ge 1} \bdelta_{\wt{\bsymb{\zeta}}^{(l)}} = \sum_{l \ge 1} \bdelta_{(\wt{\zeta}^{(l)}_1, \wt{\zeta}^{(l)}_2)}$ be a Poisson random measure on $\bbr_0^2$ with intensity measure $\bsymb{\tau} \circ \wt{\rm T}_\phi^\inv$ and independent of the random variable $W$ and the environment sequence ${\bf Y}'$. Given the sequence ${\bf Y}'$, consider a collection $({\cal R}^{(l)}_{**} : l \ge 1)$ of independent copies of the random variable ${\cal R}_{**}$ such that
\aln{
& \prob({\cal R}_{**} = r | {\bf Y}') = [\pi_{r + 1}({\bf Y}')]^\inv/ {\rm C}_5({\bf Y}') \mbox{ for all } r \ge 0 \nonumber \\
&  \mbox{ where } \pi_0({\bf Y}') = 1 \mbox{ and } \pi_j({\bf Y}') = \prod_{i = 0}^{j - 1} Y'_i \mbox { for all } j \ge 1 \mbox{ and } {\rm C}_5({\bf Y}') = \sum_{r = 1}^\infty [\pi_r({\bf Y}')]^\inv.
}
Then, it can be shown that
\aln{
{\bf N}_*^{{\tiny{\rm (lin)}}} \eqd \sum_{l = 1}^\infty [\pi_{{\cal R}_{**}^{(l)}}({\bf Y}')] \sum_{k = 1}^{Y'_{{\cal R}_{**}^{(l)}}} \bdelta_{[W {\rm C}_5({\bf Y}')]^{1/\alpha} \wt{\zeta}^{(l)}_k}. \label{eq_series_rep_weak_lim_lin_dependence}
}

\noindent{\bf Asymptotic distribution of the rightmost position.} The explicit description of the law of the rightmost position $M_{**}$ here turns out to be tricky to obtain. The first step would be using Theorem~\ref{thm_maxpos_general}. One may compute ${\bsymb \tau} \circ \wt{\rm T}_\phi^\inv ({\cal H}^{(v)}_{i_1, i_2, \ldots, i_v})$ analytically to get it. It turns out to be difficult to get a closed form from the analytic expression. But there is a probabilistic way to obtain the exact description of the asymptotic distribution of the rightmost position from an upper and lower bound involving the series representation obtained in \eqref{eq_series_rep_weak_lim_lin_dependence}. It is immediate to check that
\aln{
& \prob \big( \sum_{l = 1}^\infty \bdelta_{[W {\rm C}_5({\bf Y}')]^{1/\alpha}. \wt{\bsymb \zeta}_l} ( x, \infty) \times (x, \infty)  = 0\big) \le  \prob( {\bf N}_*(x, \infty) = 0) = \prob(M_{**} \le x) \nonumber \\
& \hspace{1cm} \le \prob \big( \sum_{l = 1}^\infty \bdelta_{[W, {\rm C}_5({\bf Y}')]^{1/\alpha} \wt{\zeta}^{(l)}_1} (x, \infty) =0 \big) . \label{eq_ineq_lin_dep_rightmost_position}
}
As $(\wt{\bsymb \zeta}^{(l)}  : l \ge 1)$ and $(\wt{\zeta}^{(l)}_1 : l \ge 1 )$ are Poisson atoms, we can condition on $W$ and ${\bf Y}'$ to compute the probabilities in the right and left hand side of the inequality \eqref{eq_ineq_lin_dep_rightmost_position}. The left hand side equals $\exptn( \exp\{ - W {\rm C}_5({\bf Y}') x^{- \alpha} \bsymb{\tau} \circ \wt{\rm T}_\phi^\inv [(1, \infty) \times (1, \infty])\})$ and the right hand side equals $\exptn( \exp \{ - W {\rm C}_5({\bf Y}')$ $x^{- \alpha} {\bsymb \tau} \circ \wt{\rm T}_\phi^\inv[(1, \infty) \times [0, \infty)] \})$. Using the description of $\bsymb{\tau}$ and $\wt{\rm T}_\phi$, one can check that $\bsymb{\tau} \circ \wt{\rm T}_\phi^\inv [(1, \infty) \times (1, \infty)] = [\max(1, \phi^\inv)]^{- \alpha}$ and $\bsymb{\tau} \circ \wt{\rm T}_\phi^\inv [(1, \infty) \times [0, \infty)]  = \phi^{\alpha}$. This computation shows that the upper and the lower bound agree if $\phi \le 1$. So we have an explicit description of the law of the rightmost position when $\phi \le 1$. If $\phi > 1$, then we do not have the explicit description for the distribution of $M_*^{\tiny{({\rm lin})}}$.

\subsection{Genealogical structure is given by a branching process in time-inhomogeneous environment} \label{subsec_bpie}


We will first consider a branching process in a time-inhomogeneous (deterministic) environment. It is well known in the literature (see the references in \cite{biggins:dsuza:1992}) that the asymptotic shape of the genealogical tree induced by BPIE is very different from that of a Galton-Watson tree or branching process in i.i.d. environment. In \cite{biggins:dsuza:1992}, the necessary conditions for the analog of Kesten-Stigum's theorem for BPIE have been derived.  Consider a collection $({\bf y}_n : n \ge 0)$ of p.m.f.s such that ${\bf y}_n = (y_n(i) : i \ge 0)$, $y_n(i) \ge 0$ and $\sum_i y_n(i) = 1$ for every $n \ge 1$. A branching process ${\bf Z} = (Z_n : n \ge 0)$ is called BPIE if for every $n \ge 1$, the size of the $n$-th generation satisfies the following distributional identity
\aln{
Z_n = \sum_{i = 1}^{Z_{n - 1}} B_{n, i},
}
where $(B_{n, i} : i \ge 1)$ is a collection of i.i.d. random variables with common p.m.f ${\bf y}_{n- 1}$ and independent of $(Z_i : 0 \le i \le n- 1)$. Let $ \mu_{n - 1} = \exptn(B_{n, 1})$ for every $n \ge 1$ and $m_n = \exptn(Z_n) = \prod_{i = 0}^{n - 1} \mu_i$ where $Z_0 = 1$ almost surely. It is clear that $(m_n^\inv Z_n : n \ge 0)$ is a non-negative martingale and hence, converges to a random variable $W$ almost surely. The main issue is to find the conditions on the environment sequence so that the limit $W$ is non-zero and non-degenerate.

\begin{fact}[Kesten-Stigum theorem for supercritical BPIE] \label{fact_bbpie_biggins_dsuza}
Let $(Z_n : n \ge 0)$ be a BPIE with $Z_0 = 1$ in the environment $({\bf y}_n : n \ge 0)$. Let the sequence $(\mu_n : n \ge 0)$ of progeny means satisfy
\aln{
\prod_{j = i}^{n + i -1} \mu_i \ge {\rm A} {\rm c}^n \mbox{ for some } {\rm A} > 0 \mbox{ and } {\rm c} > 1
\label{eq_uniform_supercriticality_bd_cond}}
for all $i \ge 1$ and $n \ge 1$. Suppose further that there exists a random variable $X$ which dominates the sequence $(\mu_{n - 1}^\inv B_{n, 1} : n \ge 1)$ and $\exptn( X \log_+ X) < \infty$. Then the sequence $(m_n^\inv Z_n : n \ge 0)$ converges almost surely to a random variable $W$ such that $\exptn(W) = 1$ and $\{ W = 0 \} = \{ \lim_{n \to \infty} Z_n = \infty\}$ almost surely.
\end{fact}
Let us now consider a BRW with regularly varying displacements and the genealogical structure given by a BPIE satisfying the assumptions in Fact~\ref{fact_bbpie_biggins_dsuza}. The uniform supercriticality condition \eqref{eq_uniform_supercriticality_bd_cond} implies almost sure convergence of $m_n^\inv \sum_{i = 1}^n Z_i$. Assuming the displacements to be i.i.d. one may employ the principle of a single large displacement to get a sequence $(b_n : n \ge 1)$ such that $\lim_{n \to \infty} m_n \prob(|X| > b_n) = 1$. Hence, the heuristics based on the ``principle of single large displacement'' holds. It is tempting to believe that the appropriately scaled rightmost position and the sequence of point processes converge weakly. We will now provide a counter-example based on a branching process in a periodic environment (BPPE).

 Let us assume for all $k \ge 0$,
\aln{
& y_{2k + 1}(2) = 1 \mbox{ and } y_{2k + 1}(j) = 0 \mbox{ for all } j \in \bbn_0 \setminus \{2\};  \nonumber \\
\mbox{ and } &  y_{2k}(1) = 1 \mbox{ and } y_{2k}(j ) = 0 \mbox{ for all } j \in \bbn_0 \setminus \{ 1\}.  \label{ass_simple_bppe}
}
A branching process in the environment $({\bf y}_k : k \ge 0)$ is the simplest example of a supercritical BPPE with period 2. It is immediate that $Z_{2k}  = Z_{2k + 1} = 2^{k}$. We define a BRW with the displacements are independent copies of the random variable $X_+$ such that $X_+ \in {\rm RV}_{- \alpha}((0, \infty), \nu_\alpha^+)$.

We denote
\aln{
{\bf G}_n = \sum_{|\sfv| = n} \bdelta_{2^{- n/(2\alpha)} S(\sfv)} \mbox{ where } 2^{n/2} \prob(2^{- n/(2\alpha)} X_+ \in \cdot) \hlconv \nu_\alpha^+(\cdot).
}
The following theorem shows that the even subsequence $({\bf G}_{2n} : n \ge 1)$ and odd subsequence $({\bf G}_{2n + 1} : n \ge 1)$ converge weakly to different limits. Recently in \cite{kriechbaum:2021} and \cite{xavier:2022}, the author proved tightness but existence of the subsequential weak limits are not yet known in this framework. The following theorem identifies the convergent subsequences and explicit description of the distribution of the weak limit.

\begin{thm} \label{thm_tighness_pp_simple_BPPE}
Under the assumptions ~\eqref{ass_simple_bppe} and $X_+ \in {\rm RV}_{- \alpha}((0, \infty), \nu_\alpha^+ )$, ${\bf G}_{2n + 1} \Rightarrow {\bf G}_*$ and ${\bf G}_{2n} \Rightarrow {\bf G}_{**}$. Furthermore,   for every $f \in C_c^+((0, \infty])$, $\exptn( \exp \{ - \int f \dtv {\bf G}_{**}\})$ equals
\alns{
&  \exp \Big\{ -  \sum_{i = 0}^\infty 2^{-i + 1 } \int_{[0, \infty)^2 \setminus \{{\bf 0}_2\}} \lambda_{**}(\dtv {\bf x})  \Big( 1 - \exp \{ - 2^i \sum_{m = 1}^2 f(x_m)\} \Big) \Big\}  
}
with $\lambda_{**}\big[(a, \infty) \times (b, \infty) \big] = [\max(a, b)]^{- \alpha} + b^{- \alpha} \mbbo_{\{a = 0\}} + a^{- \alpha} \mbbo_{\{ b = 0\}}$ for all  $(a,b) \in [0, \infty)^2$  such that  $\max(a, b) > 0$.  For every $f \in C_c^+(0, \infty])$, $\exptn(\exp \{ - \int f \dtv {\bf G}_* \})$ equals
\aln{
 &  \exp \Big\{ - \sum_{i = 0}^\infty 2^{- i + 3/2 } \int_{[0, \infty)^2 \setminus \{{\bf 0}_2\}} \lambda_*(\dtv {\bf x}) \Big( 1 - \exp \Big\{ - 2^i \sum_{m = 1}^2 f(x_m) \Big\} \Big) \Big\}  \label{eq_Laplace_G_s}\\
& = \exp \Big\{ - \sum_{i = 0}^\infty 2^{ - i + 5/2} \int_0^\infty \nu_\alpha^+(\dtv x) \big( 1 - e^{-2^i f(x)} \big) \Big\} \label{eq_Laplace_G_s_simple}
}
 with $ \lambda_*[(a, \infty) \times (b, \infty)]  = a^{- \alpha} \mbbo_{\{ b = 0\}} +  b^{- \alpha} \mbbo_{\{ a = 0\}}
 =  ( \nu_\alpha^+ \otimes \bdelta_0 + \bdelta_0 \otimes \nu_\alpha^+)[(a, \infty) \times (b, \infty)] $ for all $(a, b) \in [0, \infty)^2 \setminus \{ {\bf 0}_2\}$.  Hence, $({\bf G}_n : n \ge 1)$ is tight but does not converge weakly.
\end{thm}


The next natural question would be if we can write down the point processes ${\bf G}_*$ and ${\bf G}_{**}$ in terms of Poisson random measures. This is important to obtain an SScDPPP representation for the limit point processes and asymptotic distribution of the of $M_{2n + 1} = \sup_{|\sfv| = 2n + 1} S(\sfv)$ and $M_{2n} = \sup_{|\sfv| = 2n} S(\sfv)$ for every $n \ge 1$. It follows immediately from Theorem~\ref{thm_tighness_pp_simple_BPPE} that there exist two random variables $M_*$ and $M_{**}$ such that $2^{- (2n + 1)/(2 \alpha)}M_{2n + 1} \Rightarrow M_*$ and $2^{-n/\alpha}M_{2n} \Rightarrow M_{**}$. The interesting question is if the distribution function of $M_*$ and $M_{**}$ can be written explicitly and a representation of the limit measures using Poisson random measures plays an important role here. Consider a Poisson random measure $\sum_{l = 1}^\infty \bdelta_{E_l}$ on $(0, \infty)$ with mean measure $2^{5/2}\nu_\alpha^+$ and another Poisson random measure $\sum_{l \ge 1} \bdelta_{{\bf F}_l} = \sum_{l \ge 1} \bdelta_{(F_l^{(1)}, F_l^{(2)})}$ on $[0, \infty)^2 \setminus \{ {\bf 0}_2\}$ with mean measure $2 \lambda_{**}$. Consider a collection $(\eta_i : i \ge 1)$ of independent copies of the geometric random variable $\eta$ with probability mass function $(2^{- i - 1} : i \ge 0)$. We further assume that $(\eta_i : i \ge 1)$ is independent of the Poisson random measures $\sum_{l \ge 1} \bdelta_{E_l}$ and $\sum_{l \ge 1} \bdelta_{{\bf F}_l}$. With this notation, we are ready to provide explicit description of the weak limits ${\bf G}_*$ and ${\bf G}_{**}$ as SScDPPs.

\begin{thm}[Tightness and subsequential limits for the rightmost positions and representation of the limit random measures] \label{thm_tightness_rightmost_position}
Under the assumptions of Theorem~\ref{thm_tighness_pp_simple_BPPE}, it follows that
\aln{
{\bf G}_* \eqd \sum_{l = 1}^\infty 2^{\eta_l} \bdelta_{ E_l} \mbox{ and } {\bf G}_{**} \eqd \sum_{l = 1}^\infty 2^{\eta_l}(\bdelta_{F_l^{(1)}} + \bdelta_{F_l^{(2)}}).
}
Furthermore, for all $x > 0$, we have
\aln{
& \prob(M_* \le x)  =  e^{- 2^{5/2} x^{- \alpha}} \mbox{ and } \prob(M_{**} \le x) =  e^{- 6 x^{- \alpha}}.
}
Hence, the sequence $(2^{- n/(2\alpha)} M_n : n \ge 1 )$ is tight but does not converge weakly.
\end{thm}

\begin{remark}[No weak convergence for BPPE and hence, BPIE]
Theorem~\ref{thm_tighness_pp_simple_BPPE} indicates that a version of Theorem~\ref{thm:main:thm} does not hold where the genealogical structure is given by a BPPE. But, in our simple example, the normalized extreme positions are tight. We believe that the sequence of the point processes remains tight while the genealogical tree is a BPIE. But the rigorous proof of this belief may turn out to be a challenging task. It is clear from Theorem~\ref{thm_tighness_pp_simple_BPPE} that it is not possible to obtain weak convergence of the point processes for BPIE. After proving tightness and assuming the genealogical tree to be BPPE, it seems to be an interesting and challenging question to identify the subsequential limits.

\end{remark}


\begin{remark}[Branching process in a Markovian environment] \label{remark_BPME_cojec}
In the previous remark, we have seen that Theorem~\ref{thm:main:thm} can not be extended when the genealogical tree is a branching process in a time-varying (deterministic) environment. This observation does not exclude the possibility of extension when the genealogical tree is a branching process in a random environment (other than an i.i.d. environment).  The branching process in an ergodic reversible Markovian environment is of special interest as the exchangeable environment can be dissected into i.i.d. regeneration cycles. Suppose that one is interested in the scaling/normalizing sequence which depends only upon the genealogical tree. Then it is tempting to believe that there exists a (random) sequence $(C_n : n \ge 1)$ such that $({\bf N}_n = \sum_{|\sfv| = n} \bdelta_{C_n^\inv S(\sfv)} : n \ge 1)$ turns out to be tight. But for any such sequence $(C_n : n \ge 1)$, ${\bf N}_n$ may not converge weakly. As we expect different phenomena, in this case, Theorem~\ref{thm:main:thm} can not be extended directly in this case. Nevertheless, we believe that the main theorem in this article will turn out to be crucial in the formalization of this phenomenon.
\end{remark}

\section{Sketch of the proof of Theorem~\ref{thm:main:thm} and proofs of Theorems~\ref{thm_tighness_pp_simple_BPPE} and \ref{thm_tightness_rightmost_position}} \label{sec_proof_main_thm}

The theorem follows if one can show that limit of the Laplace functional of ${\bf N}_n$ exists as $n \to \infty$ and the limit of the Laplace functional is a Laplace functional of ${\bf N}_*$.
The derivation is long, tedious, and similar to the proof of Theorem~2.6 in \cite{bhattacharya:hazra:roy:2017}. In the first three steps, we shall introduce three new point processes (which involve the displacements rather than the positions) to approximate ${\bf N}_n$. The main aim of these steps is to identify the location (in the genealogical tree) of the largest displacements with high probability as the heuristics based on the largest displacements turn out to be correct. In the final step, the limit of the Laplace functional of ${\bf N}_n$ will be computed by counting the number of the positions affected by the large displacements in the $n$-th generation. We shall now briefly mention the steps for the sake of completeness.

The first step would be formalization of the fact that at the most one displacement on a geodesic path $o \mapsto \sfv$ can be large. Informally, ${\bf N}_n = \sum_{|\sfv| = n} \bdelta_{B_n^\inv \sum_{\sfu \in o \mapsto \sfv} X(\sfu)} \approx \sum_{|\sfv| = n} \sum_{\sfu \in o \mapsto \sfv} \bdelta_{B_n^\inv X(\sfu)}$ where `$\approx$' stands for `well approximated by'. The precise statement is as follows.

\begin{propn}\label{propn_old_geodesic_path}
If  the assumptions of Theorem~\ref{thm:main:thm} hold, then for every $\epsilon > 0$, $\limsup_{n \to \infty}$  $\prob^* [\dvague({\bf N}_n, \wtbfN_n) > \epsilon] = 0$.
\end{propn}

Once we approximate the point process ${\bf N}_n$ by a point process composed of the displacements, we need location of the large displacements (of absolute magnitude $B_n$) in the genealogical tree with high probability.  As the genealogical tree has geometric growth and each bunch of displacements, coming from  the same parents has the same probability to be the large, the first $(n - \varrho)$ generations will not contain the large ones with high probability when $\rho$ is an intermediate sequence.

\begin{lemma}[Forest contains the the large displacements with high probability] \label{lemma_cutting}
If the assumptions in Theorem~\ref{thm:main:thm} hold, then we have $\lim_{\varrho \to \infty} \limsup_{n \to \infty} \prob^* \big[ \dvague(\wt{\bf N}_n , \wt{\bf N}_{n, \varrho}) > \epsilon \big] = 0$ for all $\epsilon > 0$ where
\[ \wt{\bf N}_{n, \varrho}= \sum_{|\uv| = n} ~~ \sum_{\substack{\sfu \in o \mapsto \sfv\\ |\sfu \mapsto \sfv| \le \varrho - 1}} \bdelta_{B_n^\inv X(\sfu)}. \]
\end{lemma}

After cutting the tree into forest, it is clear that each of the subtrees (with root at the $(n - \varrho)$-th generation of ${\sf T}$) can contain the large displacement with equal probability. Roughly speaking, it turns out that the number subtrees containing the large displacements approximately Poisson distributed. Here the main differences from the earlier work are that the subtrees are dependent through the environment and the law of the subtrees depend on the scaled-displacements through scaling. It is believable that these issues can be taken care of by conditioning on the environment. Focusing on a single subtree, we can see that each particle can produce a bunch of large displacements with equal probability. We want to show that the particles with exceptionally large number of children can not produce the large displacements with high probability. To formalize this, we need some new notation. Let $\wh{T}_i$ denotes the $i$-th subtree in the forest with root ${\sf r}_i$ for all $i \ge 1$. If $\sfu \in \wh{\sf T}_i$, then $A(\sfu)$ denotes the number of descendants of $\sfu$ in the $\varrho$-th generation of $\wh{T}_i$. With these notations, it is easy to see that $\wt{\bf N}_{n, \varrho} = \sum_{i = 1}^{Z_{n - \varrho}} \sum_{\sfu \in \wh{T}_i \setminus \{ {\sf r}_i\}} \bdelta_{B_n^\inv X(\sfu)}$. As the BPRE is supercritical, there exists a large integer $\vartheta$ such that each subtree remain supercritical even though each particle is allowed to produce at the most $\vartheta$. We call $\vartheta$ a pruning threshold. We modify each subtree such that each particle in the subtree has at the most $\vartheta$ children. After modification, $\wh{T}_i$ is called pruned subtree and denoted by $\ol{\sf T}_i$.

\begin{propn}\label{propn_pruning}
If the assumptions in Theorem~\ref{thm:main:thm} hold, then for all $\varrho \ge 1$ and $\epsilon > 0$, $\lim_{\vartheta \to \infty}$ $  \limsup_{n \to \infty}$  $\prob^* [ \dvague(\wt{\bf N}_{n, \varrho, \vartheta}, \wt{\bf N}_{n,\varrho }) > \epsilon] = 0$ where
\alns{
\wt{\bf N}_{n, \varrho, \vartheta} = \sum_{i = 1}^{Z_{n - \varrho}} \sum_{\sfu \in \ol{\sf T}_i \setminus \{{\sf r}_i\}} A^{(\vartheta)} (\sfu) \bdelta_{B_n^\inv X(\sfu)}
}
and $A^{(\vartheta)}(\sfu)$ denotes the descendants of the vertex $\sfu$ in the $\varrho$-th generation of $\ol{\sf T}_i$ if $\sfu \in \ol{\sf T}_i$.
\end{propn}

\noindent{{\bf Step~4}. (Regularization of the forest).} In this step, we add new children to each of the vertex, so that it has exactly $\vartheta$ children in the next generation. The modified subtree is a $\vartheta$-regular tree and denoted by $\wt{\sf T}^{(i)}$ for every $i \ge 1$.  If $\sfu \in \wt{\sf T}^{(i)} \setminus \ol{\sf T}^{(i)}$, then $A^{(\vartheta)}(\sfu) := 0$. For every $\sfu  \in \wt{\sf T}^{(i)}$, we replace the displacements of its $\vartheta$ children by an independent copy of the random vector $(X_1, X_2, \ldots, X_\vartheta)$. The displacement attached to the vertex $\sfu$ will be denoted by $X'(\sfu)$.  It is clear that
\aln{
\wtbfN_{n, \varrho, \vartheta} \eqd \sum_{i = 1}^{|{\sf D}_{n - \varrho}|} \sum_{\sfu \in \wt{\sf T}^{(i)} \setminus \{{\sf r}_i\}} A^{(\vartheta)}(\sfu) \bdelta_{B_n^\inv X'(\sfu)}. \label{eq_pointprocess_regularized}
}
With a slight abuse on notation, we use $\wtbfN_{n, \varrho, \vartheta}$ to denote the point process in the right hand side of \eqref{eq_pointprocess_regularized}. In the following proposition we show that $\wtbfN_{n, \varrho, \vartheta}$ converges weakly to the desired limit.

\begin{propn} \label{propn:laplace:regular:point:process}
Recall ${\bf N}_*$ from \eqref{characterization}. Under the assumptions of Theorem~\ref{thm:main:thm}, for every $f \in \ccr$,
\aln{
\lim_{\varrho \to \infty} \lim_{\vartheta \to \infty} \lim_{n \to \infty} \exptn^* \big[ \exp \big\{ - \int f \dtv \wtbfN_{n, \varrho, \vartheta} \big\}  \big]  = \exptn^* \big[ \exp \big\{ - \int f \dtv {\bf N}_* \big\} \big].
}
\end{propn}

The detailed proofs of these propositions are included in the Appendix.

\subsection{Proof of Theorem~\ref{thm_tighness_pp_simple_BPPE}}

It is clear that we need to derive the weak limits of the following two sequences $({\bf G}_{2n + 1} : n \ge 0)$ and $({\bf G}_{2n} : n \ge 0)$.
In Figure 1 below, we have used blue and red for the vertices in the odd and even generations, respectively. We have mentioned the displacements attached to the vertices in the picture for the first few generations.

To understand the asymptotic behavior of the point process ${\bf G}_{2n}$, one needs to understand the asymptotic behavior of the positions of the red particles. Note that the genealogical structure of the red particles is also given by a binary tree $\bbt^{(r)}$. But the progeny point process ${\mathscr L}_{**} \eqd \bdelta_{X(1) + X(1, 1)} + \bdelta_{X(1) + X(1, 2)}$ has dependent atoms. This fits the framework of \cite{bhattacharya:hazra:roy:2017}. Here, we have to compute the $\mathbb{M}_0$-limit of the sequence $(2^{n + 1} \prob[ 2^{ - n/\alpha} (X(1) + X(1, 1), X(1) + X(1, 2)) \in \cdot] : ~ n \ge 1)$ which is difficult to guess due to the dependence. We shall first note that $ 2^{n + 1} \prob \big[ 2^{- n/\alpha} \big( X_+ (1), X_+(1, 1), X_+(1, 2) \big) \in \cdot  \big] \hlconv 2 \boldsymbol{\tau}_{**}(\cdot)$ where $\boldsymbol{\tau}_{**} = \nu_\alpha^+ \otimes \bdelta_0 \otimes \bdelta_0 + \bdelta_0 \otimes \nu_\alpha^+ \otimes \bdelta_0 + \bdelta_0 \otimes \bdelta_0 \otimes \nu_\alpha^+$. We can apply the mapping theorem (Theorem~2.3 in \cite{lindskog:resnick:roy:2014}) for $\mathbb{M}_0$-convergence to derive
\aln{
& 2^{n + 1} \prob \big[ 2^{n /\alpha} \big( X_+ (1) + X_+(1, 1), X_+(1) + X_+(1, 2) \big) \in \cdot  \big]  \hlconv 2 \lambda_{**}(\cdot)\nonumber \\
&  \mbox{where } \lambda_{**}(\cdot) = \big[ \nu_\alpha^+ (\{ x > 0: (x, x) \in \cdot \}) + \nu_\alpha^+ \otimes \bdelta_0 (\cdot) + \bdelta_0 \otimes \nu_\alpha^+ (\cdot) \big].  \label{eq_defn_lambda_ss_proof}
}
The rest of the computation follows from the expression in (4.44) in \cite{bhattacharya:hazra:roy:2017} replacing $\lambda$ by $2 \lambda_{**}$ and the other necessary modifications.

For the asymptotic behavior of ${\bf G}_{2n + 1}$, one needs to understand the positions of the blue particles. It is convincing that with high probability, $X_+(1)$ is not large enough when $n$ is large enough. So we can ignore the contribution of the displacement $X_+(1)$. If we consider only the blue particles, then we can see that they form a binary tree $\bbt^{(b)}$. So we have a BRW (modified) indexed by a binary tree where the progeny point process ${\mathscr L}_* \eqd \bdelta_{X_+(1, 1) + X_+(1, 1, 1)} + \bdelta_{X_+(1,2) + X_+(1, 2, 1)}$. It is clear from the assumptions that $X_+(1,1) + X_+(1, 1, 1) $ and $X_+(1,2) + X_+(1, 2, 1)$ are i.i.d. Hence, we can use Theorem~ in \cite{bhattacharya:hazra:roy:2016} to obtain the weak limit of ${\bf G}_{2n + 1}$.

\noindent \begin{minipage}{0.6 \textwidth}
It is clear that the $(2n + 1)$-th generation of the original tree contains $2^n$ blue particles and corresponds to the $n$-th generation of $\bbt^{(b)}$ and implies that the size of $\bbt^{(b)}$ is of order $2^{n + 1}$ upto the $n$-th generation.  We shall be done if we can compute the $\mathbb{M}_0$-limit of the sequence of measures $( 2^{n + 1} \prob[ 2^{- (2n + 1)/(2\alpha)} (X(1, 1) + X(1, 1, 1), X(1, 2) + X(1, 2, 1)) \in \cdot ]: ~ n \ge 1)$.  As $X(1,1) + X(1, 1, 1) $ and $X(1,2) + X(1, 2, 1)$ are i.i.d. it is easy to guess that $\mathbb{M}_0$ limit of the sequence of measures should concentrate on the axes. But we also need the algebraic factor in front of the limit measure for an explicit description which needs some careful computation.  Note that (due to i.i.d. structure)
\end{minipage}
\begin{minipage}{0.4 \textwidth}
\hspace{4cm}
\begin{tikzpicture}[transform canvas={scale=.7}]

\draw[red](-1,3) node{\textbullet};

\draw[-](-1, 3) -- (-1,2);
\draw[blue](-1,2) node{\textbullet};
\draw[](-1,2) node[anchor= north]{$X_+(1)$};

\draw[-](-1,2) -- (-3, 1);
\draw[-](-1,2) -- (1, 1);
\draw[red](-3, 1) node{\textbullet};
\draw(-3, 1) node[anchor = east]{\text{\small $X_+(1,1)$}};
\draw[red](1,1) node{\textbullet};
\draw (1,1) node[anchor = west]{\text{ \small $X_+(1,2)$}};

\draw[-](-3, 1) -- (-3, 0);
\draw[blue] (-3,0) node{\textbullet};
\draw(-3, 0) node[anchor = east]{\text{\small $X_+(1,1,1)$}};

\draw[-](1, 1) -- (1, 0);
\draw[blue](1,0) node{\textbullet};
\draw(1, 0) node[anchor = east]{\text{\small $X_+(1,2,1)$}};


\draw[-](-3,0) -- (-3.5,-1);
\draw[-](-3,0) -- (-2.5, - 1);
\draw[red](-3.5,-1) node{\textbullet};
\draw[red](-2.5, -1) node{\textbullet};

\draw[-](1,0) -- (1.5, -1);
\draw[-](1,0) -- (.5,-1);
\draw[red](1.5, -1) node{\textbullet};
\draw[red](.5, -1) node{\textbullet};


\draw[-] (-3.5, - 1) -- (-3.5, - 2);
\draw[blue](-3.5, - 2) node{\textbullet};
\draw[-](-2.5, - 1) -- (-2.5, - 2);
\draw[blue] (-2.5, - 2) node{\textbullet};

\draw[-](1.5, - 1) -- (1.5, - 2);
\draw[blue](1.5, -2) node{\textbullet};
\draw[-](.5, - 1) -- (.5, - 2);
\draw[blue](.5, -2) node{\textbullet};

\draw(-1,-3) node{\text{Figure~1: Genealogical tree for the BPPE.}};
 \end{tikzpicture}
 \end{minipage}

\alns{
& 2^{n + 1} \prob \big[ 2^{- (2n + 1)/(2 \alpha)} \big( ~X_+(1, 1), X_+(1, 1, 1), X_+(1, 2), X_+(1, 2, 1) ~ \big) \in \cdot \big] \hlconv  2^{1/2} \boldsymbol{\tau}_*( \cdot) \\
& \mbox{ where }  \boldsymbol{\tau}_* = \nu_\alpha^+ \otimes \bdelta_0 \otimes \bdelta_0 \otimes \bdelta_0 + \bdelta_0 \otimes \nu_\alpha^+ \otimes \bdelta_0 \otimes \bdelta_0 + \bdelta_0 \otimes \bdelta_0 \otimes \nu_\alpha^+ \otimes \bdelta_0 + \bdelta_0 \otimes \bdelta_0 \otimes \bdelta_0 \otimes \nu_\alpha^+.
}
Then one can use the mapping theorem (Theorem~2.3 in \cite{lindskog:resnick:roy:2014})  for $\mathbb{M}_0$ convergence and some algebra to see that $ 2^{n + 1} \prob[ 2^{- (2n + 1)/(2\alpha)} (X(1, 1) + X(1, 1, 1), X(1, 2) + X(1, 2, 1)) \in \cdot ] \hlconv \lambda_* (\cdot)$ where $ \lambda_*(\cdot) = 2^{3/2} \big( \nu_\alpha^+ \otimes \bdelta_0 + \bdelta_0 \otimes \nu_\alpha^+ \big)(\cdot)$.
The rest of the computation follows from the expression of the Laplace functional obtained in (2.6) in \cite{bhattacharya:hazra:roy:2016}  after appropriate modifications.

\subsection{Proof of Theorem~\ref{thm_tightness_rightmost_position}}

The series representation of ${\bf G}_*$ follows immediately from the remark~2.4 in \cite{bhattacharya:hazra:roy:2016} and ${\bf G}_{**}$ follows from Theorem~ 2.6 in \cite{bhattacharya:hazra:roy:2017} (see also Subsection~4.4 for a detailed derivation). It is clear that $\prob(M_* \le x) = \prob( {\bf G}_* (x, \infty) = 0) = \exp\{ - 2^{5/2} \nu_\alpha^+(x, \infty)\}$ as $\{ {\bf G}_*(x, \infty) = 0\} = \{ \sum_{l = 1}^\infty \bdelta_{E_l} (x, \infty) = 0\}$.   Following the same arguments, we can see that
\alns{
\prob \big( M_{**} \le x \big) & = \prob \Big( \sum_{l = 1}^\infty (\bdelta_{F_l^{(1)}} + \bdelta_{F_l^{(2)}}) (x, \infty) = 0 \Big)  = \prob \Big( \sum_{l = 1}^\infty \bdelta_{(F_l^{(1)}, F_l^{(2)})} \big( [0, \infty)^2 \setminus [0, x]^2 \big) = 0 \Big) \\
& = \exp \big\{ - 2 \lambda_{**}  \big( (x, \infty) \times [0, \infty) \cup [0, \infty) \times (x, \infty) \big) \big\}.
}
One can use the description of the measure $\lambda_{**}$ in \eqref{eq_defn_lambda_ss_proof} to get an explicit expression for the exponent.

\section*{Acknowledgement}
This research is partially supported by
Polish National Science Centre Grant No. 2018/29/B/ST1/00756.

\bibliographystyle{abbrvnat}

\section{Appendix}
We shall provide detailed proofs of the auxiliary results for Theorem~\ref{thm:main:thm}, mentioned in Section~\ref{sec_proof_main_thm}.

\section{Notations and facts}

Recall that ${\bf Y} = (Y_i : i \ge 0)$ denotes the environment sequence and $\exptn(Z_n | {\bf Y}) = \pi_n$ for all $n \ge 1$. As the environment is i.i.d. one can use SLLN to have
\aln{
n^\inv \log \pi_n \asconv \exptn(\log \exptn(Z_1|Y_0)) > 0. \label{eq_geometric_growth_BPRE}
}
We assume that the displacements coming from the same parent is jointly regularly varying on the space $\bbr^{\bbn}_{{\bf 0}}$. See Assumption~2.3 in page~13 of the aforementioned reference. We defined
\aln{
B_n := \inf \{ s > 0 : \prob( |X_1| > s) < \pi_n^\inv \}.
\label{eq_defn_Bn}
}
It has been shown that
\aln{
& \pi_n \prob \Big( B_n^\inv {\bf X} \in \cdot \Big| {\bf Y}_{0~:~n - 1} \Big) \hlconv \bsymb{\nu}_*(\cdot) \mbox{ almost surely }   \label{eq_mzero_conv_crucial}\\
& \mbox{ where } \bsymb{\nu}_* = \bsymb{\nu}( \cdot ) /\bsymb{\nu} ({\sf E}_1) \mbox{ and } {\sf E}_1 = \big( \big\{ (x_i : i \ge 1) : |x_1| > 1 \big\} \big). \nonumber
}

\medskip

\noindent{\bf A fact from the theory of $\mathbb{M}_0$ convergence.} In order to derive the limit of the Laplace functional after regularization, the connection between the finite-dimensional convergence (a consequence of Theorem~4.1 in \cite{lindskog:resnick:roy:2014}) and $\mbbM_0$ convergence on $\bbrz$ will be useful for us.
\begin{fact} \label{fact_finite_dimensional_conv_mzero_conv}
Let ${\bf 0}_d$ denote the origin of $\bbr^d$ for all $d \ge 1$. Consider the projection map $\proj_d: \bbr^\bbn \to \bbr^d$ defined by $\proj_d ((x_i : i \ge 1)) = (x_i : 1 \le i \le d)$ for every $d \ge 1$. Then $\lambda_n \hlconv \lambda$ in $\mbbM(\bbrz)$ if and only if for all $d \ge 1$, we have
\aln{
\lambda_n \circ \proj_d^\inv \hlconv \lambda \circ \proj_d^\inv \mbox{ in } \mbbM(\bbr_{{\bf 0}}^d) \mbox{ where } \bbr_{{\bf 0}}^d = \bbr^d \setminus \{{\bf 0}_d\}.
}
\end{fact}

\section{Proofs of the auxiliary results in Section~3}

Here, we shall briefly discuss the proofs of the auxiliary results. Some of the auxiliary results needed for these proofs will be proved later in Section~\ref{sec_further_proofs}.

\subsection{Proof of Proposition~3.1}
 Consider $f \in \ccr$ such that ${\rm supp}(f) := \{x \in \ol{\bbr}: f(x) \neq 0\} \subset \{ x \in \bbr: |x| > \vep\}$ for some $\vep > 0$. The proposition follows if one can show that for all $\ol{\epsilon} > 0$,
\aln{
\limsup_{n \to \infty} \prob^*[ {\sf E}_{n , \ol{\epsilon}}] = 0 \mbox{ where } {\sf E}_{n , \ol{\epsilon}} = \big\{~ \big| \sum_{|\sfv| = n} f(B_n^\inv S(\sfv)) - \sum_{|\sfv| = n} \sum_{\sfu \in o \mapsto \sfv} f(B_n^\inv X(\sfu)~) ~\big| > \ol{\epsilon} \big\}. \label{eq_aim_old_geodesic_path}
}
The key to the proof is to show that at the most one large (of the order of $B_n$) displacement can occur on a geodesic path $o \mapsto \sfv$. To be mathematically precise, we define
\aln{
{\sf G}_{n, \eta} = \Big\{ \bigcup_{|\sfv| = n} \Big\{ \sum_{\sfu \in o \mapsto \sfv} \bdelta_{B_n^\inv |X(\sfu)|} (\eta/n, \infty) \ge 2 \Big\}~~ \Big\}^c \label{eq_defn_one_large_disp_event}
}
where ${\sf B}^c$ is the complement of the event ${\sf B}$ and $\eta$ is a positive number. One can show that $\limsup_{n \to \infty} \prob^* \big[ {\sf G}_{n, \eta} \cap {\sf E}_{n , \ol{\epsilon}} \big] = 0$  for an appropriately chosen $\eta$ molding proof of (4.9) in \cite{bhattacharya:hazra:roy:2016}. Proof of \eqref{eq_aim_old_geodesic_path} follows from
\aln{
\limsup_{n \to \infty} \prob^* [( {\sf G}_{n, \eta})^c  ] = 0 \mbox{ for all } \eta > 0. \label{eq_aim_old_event}
}
We shall now provide details for the proof of \eqref{eq_aim_old_event}. As $\prob^*$ is absolutely continuous with respect to $\prob$, we can replace $\prob^*$ in \eqref{eq_aim_old_event} by $\prob$. Observe that $\prob^*[ ({\sf G}_{n, \eta})^c]$ is bounded from above by
\aln{
\exptn \Big[ 1 \wedge \Big( \exptn _{\bf Y} \Big[ \sum_{|\sfv| = n} \prob_{\bf Y} \big[ \sum_{\sfu \in o \mapsto \sfv} \bdelta_{B_n^\inv |X(\sfu)|} (\eta/n, \infty) \ge 2~ \big| ~ Z_n ~\big] \Big] \Big) \Big] \label{eq_old_upp_bound_event}
}
using the union bound conditioned on the environment ${\bf Y}$. Further, if the genealogical tree ${\sf T}$ extincts before the $n$-th generation, then \eqref{eq_aim_old_event} holds trivially. So it is enough to deal with the conditional expectation only when $Z_n > 0$. On a typical geodesic path $o \mapsto \sfv$, the displacements are i.i.d. One can use this fact along with the Markov inequality to obtain the following upper bound for the conditional expectation in \eqref{eq_old_upp_bound_event}
\aln{
{\rm const.} ~\pi_n \big[ \prob_{\bf Y} \big( |X_1| > B_n \eta/n \big) \big]^2 \label{eq_old_upp_bound_cond_exptn}
}
where $\pi_n$ denotes the quenched size of the $n$-th generation. Due to the randomness of $B_n$ inherited from the environment, we introduce another event
\aln{
{\sf J}_{n, \eta} (\wt{\vep}) : = \{ |\pi_n \prob(|X_1| > B_n \eta | {\bf Y}_{0~:~n-1}) - \eta^{- \alpha}| < \wt{\vep}\} \mbox{ for every } \wt{\vep} > 0.  \label{eq_old_defn_Jn_eta_wt_vep}
}
It follows from \eqref{eq_mzero_conv_crucial} (combined with the continuous mapping theorem for $\mathbb{M}_0$-convergence and Fact~\ref{fact_finite_dimensional_conv_mzero_conv}) that $\prob[ {\sf J}_{n, \eta}(\wt{\vep}) ] = o(1)$ as $n \to \infty$ for all $\eta > 0$ and $\wt{\vep} > 0$.  So, we must analyze the upper bound in \eqref{eq_old_upp_bound_cond_exptn} on the event ${\sf J}_{n, \eta}(\wt{\vep})$
and show that
\aln{
\limsup_{n \to \infty} \exptn \Big[ 1 \wedge \Big( \mbbo_{{\sf J}_{n, \eta}(\wt{\vep})} ~~\pi_n ~~[ \prob_{\bf Y}( |X_1| > B_n \eta/n )]^2 \Big) \Big] = 0. \label{eq_old_event_aim_reduced}
}
Our first task is to estimate $\prob_{\bf Y}(|X_1| > B_n \eta)$ which will be done by comparing this probability to $\prob_{\bf Y}(|X_1| > B_n \eta) $ via {\it stochastic version of the Potter's bound} (see Proposition~0.8(ii) in \cite{resnick:1987}). To invoke this bound, one needs to ensure that $B_n$ is growing faster than any linear function of $n$ with high probability.
To see this note that from \eqref{eq_defn_Bn} and Proposition~0.8(v) in \cite{resnick:1987}, it follows that $B_n = \pi_n^{1/\alpha} \wt{\rm L}(\pi_n)$ conditioned on ${\bf Y}$ where $\wt{\rm L}$ is a slowly varying function. Let $\mu = \exptn(\log [\exptn(Z_1 | Y_0)]) >0$ be the almost sure limit of $n^\inv \log \pi_n$ (see \eqref{eq_geometric_growth_BPRE}). Fix an $\epsilon_1 \in (0, e^\mu - 1)$ and define ${\sf H}_{n, \epsilon_1} = \{ |\pi_n^{1/n} - e^\mu| < \epsilon_1\}$. In view of the SLLN, it is enough to show that $n^\inv B_n$ can be arbitrarily large on the event ${\sf H}_{n, \epsilon_1}$ as $\lim_{n \to \infty } \prob( {\sf H}_{n, \epsilon_1} ) = 0$. On the event ${\sf H}_{n, \epsilon_1}$, it is immediate to check that $\pi_n > (e^\mu - 1)^{n/\alpha} \uparrow \infty$ and hence, Potter's bound (Proposition~0.8(ii) in \cite{resnick:1987}) can be used for $\wt{\rm L}$. Fix $\epsilon_2 \in (0, 1 \wedge \alpha^\inv )$  and for large enough $n$, $B_n  \ge (1 - \epsilon_2) (e^\mu - \epsilon_1)^{n(1/\alpha - \epsilon_2)}$ on the event ${\sf H}_{n, \epsilon_1}$. Hence we derived the following lemma.

\begin{lemma}\label{lemma_geometric_growth_Bn}
Under the assumptions of Theorem~2.3, $n^\inv B_n \probconv \infty$.
\end{lemma}
Fix $\vep_1 \in (0, 1)$. It follows from Proposition 0.8(ii) in \cite{resnick:1987} that for every $x >1$, there exists a $t_0$ such that $\prob(|X_1| > t) \le ( 1 - \vep_1)t^{ \alpha + \vep_1} \prob(|X_1| > x t) $ for all $t \ge t_0$.  Define $\wt{\sf K}_{n, t_0} := \{n^\inv B_n > t_0\} $. Then on the event $\wt{\sf K}_{n, t_0}$, we have  $\prob_{\bf Y}(|X_1| > B_n \eta/ n) \le (1- \vep_1)$  $n^{ \alpha + \vep_1} \prob_{\bf Y} (|X_1| > B_n \eta)$.
Similarly on the event ${\sf J}_{n, \eta}(\wt{\vep})$, we have $[\pi_n \prob_{\bf Y}(|X_1| > B_n \eta ) ]^2 \le (\eta^{- \alpha} + \wt{\vep})^2$.
Then 
the term inside the expectation in \eqref{eq_old_event_aim_reduced} can be decomposed into two terms based on the event $\wt{\sf K}_{n, t_0}$. After some algebraic adjustments and using the last two upper bounds, we are left with the following upper bound to the expectation in \eqref{eq_old_event_aim_reduced}
\aln{
& \exptn \Big[ \min \Big( 1, (\eta^{- \alpha} + \wt{\vep})^2 (1 - \vep_1)^2 (\pi_n^\inv n^{2\alpha + 2 \vep_1}) \mbbo_{\wt{\sf K}_{n, t_0}} \Big) \Big]  \nonumber \\
& \hspace{1cm} + \exptn \Big[ \min \Big(1, (\eta^{- \alpha} + \wt{\vep})^2 \Big[ \frac{\prob_{\bf Y}(|X_1| > B_n \eta/n)}{\prob_{\bf Y}(|X_1| > B_n \eta)} \Big]^2 \pi_n^\inv \mbbo_{\wt{\sf K}^c_{n, t_0}} \Big) \Big] =: {\rm T}_n^{(1)} + {\rm T}_n^{(2)}
}
for large enough $n$.  Geometric growth of $\pi_n$ implies that $n^{2 \alpha + 2 \delta} \pi_n^\inv \asconv 0$. So we can use the dominated convergence theorem to conclude ${\rm T}_n^{(1)} = o(1)$. The term inside the expectation of ${\rm T}_n^{(2)}$ is bounded from above by $\mbbo_{\wt{\sf K}^c_{n , t_0}}$ almost surely. It follows immediately from Lemma~\ref{lemma_geometric_growth_Bn} that $T_n^{(2)} = o(1)$. Hence, the claim in \eqref{eq_old_event_aim_reduced} follows.

\subsection{Proof of Lemma~3.2}
Consider a function $f \in \ccr$ such that ${\rm supp}(f) \subset \{x : |x| > \vep\}$. To prove the lemma, it is enough to show that $\lim_{\varrho \to \infty}\limsup_{n \to \infty} \prob^* (|\wtbfN_n(f) - \wtbfN_{n, \varrho}(f)| > \epsilon) =0$. The key to this proof is the following event
\aln{
{\sf O}_{n, \eta, \varrho} : = \Big\{ \sum_{|\uu| \le n - \varrho} \bdelta_{|\Xu|}(B_n \eta, \infty) \ge 1 \Big\}
}
as $\prob^* ( \{ |\wtbfN_n(f) - \wtbfN_{n, \varrho}(f)| > \epsilon \} \cap {\sf O}_{n, \eta, \varrho}^c) = 0$ if $\eta < \vep$. So it is enough to show that $\lim_{\varrho \to \infty}$ $ \lim_{n \to \infty} \prob({\sf O}_{n, \eta, \varrho}) = 0$. We can replace $\prob^*$ by $\prob$ as $\prob^*$ is absolutely continuous with respect to $\prob$.

Conditioned on the environment ${\bf Y}$, one may use the union bound to obtain the following upper bound for $\prob \big( {\sf O}_{n, \eta, \varrho} \big)$
\aln{
& \exptn \Big[ \min \Big( 1, \exptn_{\bf Y} \big[ \sum_{i = 1}^{n - \varrho} Z_i \big] \prob_{\bf Y}(|X_1| > B_n \eta) \Big) \Big]  \nonumber \\
& \le \exptn \Big[ \min \Big( 1, \big[ \pi_n \prob_{\bf Y}(|X_1| > B_n \eta)\big] \big[ \sum_{i = 1}^{n - \varrho} (\pi_i/ \pi_n) \big] \Big) \mbbo_{{\sf J}_{n, \eta}(\wt{\vep})} \Big] + \prob \big( {\sf J}_{n, \eta}^c (\wt{\vep}) \big) \nonumber \\
& \le \exptn \Big[ \min\Big( 1, (\eta^{- \alpha} + \wt{\vep}) \big[ \pi_n^\inv \sum_{i = 1}^{n - \varrho} \pi_i \big] \Big) \Big] + o(1) \label{eq_final_upper_bound_cutting_event}
}
where ${\sf J}_{n, \eta}(\wt{\vep})$ is introduced in \eqref{eq_old_defn_Jn_eta_wt_vep}. We now use the exchangeability property of the environment to deduce the following
\aln{
\sum_{i = 0}^{n - \varrho} \pi_i/\pi_n = \sum_{i = 0}^{n - \varrho} \big[ \prod_{j = i + 1}^{n - 1} \exptn(\xi | Y_j) \big]^\inv \eqd \sum_{i = \varrho}^{n - 1} \Big( \prod_{j = 0}^{i - 1} \exptn_{Y_j}(\xi) \Big)^\inv = \sum_{i = \varrho}^{n - 1} \pi_i^\inv. \label{eq_excheability_distn_identity}
}
Recall from Remark~1.3 that $\sum_{i \ge 1} \pi_i^\inv < \infty$. Therefore,  the right hand side of \eqref{eq_excheability_distn_identity} vanishes almost surely as $n \to \infty$ and $\varrho \to \infty$.  We can use now dominated convergence theorem to conclude that the first term in \eqref{eq_final_upper_bound_cutting_event} converges to zero as $n \to \infty$ and $\varrho \to \infty$.

\subsection{Proof of Proposition~3.3}
Let $f \in \ccr$ be a function such that ${\rm supp}(f) = \{x > 0: |x| > \vep\}$. Then it is enough to show that
\aln{
\lim_{\vartheta \to \infty} \limsup_{n \to \infty} \prob^* \Big( \Big| \wtbfN_{n, \varrho, \vartheta} (f) - \wtbfN_{n, \varrho}(f) \Big| > \epsilon \Big) = 0 \label{eq_aim_pruning}
}
for every $\varrho \ge 1$. To prove \eqref{eq_aim_pruning}, we shall construct pre-pruned subtrees through a marking scheme and the pre-pruned subtrees will be denoted by $({\sf T}_{\varrho, i}^\diamond: 1 \le i \le Z_{n - \varrho})$. To be more specific and precise, one can design a marking scheme (following Subsection~(4.5.1) in \cite{bhattacharya:2018a}) to separate out the vertices (in $\ol{\sf T}_i$)  to be included in $\wh{\sf T}_i$. Then we can consider an intermediate point process ${\bf N}_{n, \varrho, \vartheta}^\diamond$ indexed by  $({\sf T}_{\varrho, i}^\diamond: 1 \le i \le Z_{n - \varrho})$. The proof of \eqref{eq_aim_pruning} follows if we can show
\aln{
& \lim_{\vartheta \to \infty} \limsup_{n \to \infty} \prob^* \Big( \Big| \wtbfN_{n, \varrho, \vartheta}(f) - {\bf N}_{n, \varrho, \vartheta}^\diamond (f) \Big| > \epsilon/2 \Big) = 0 \label{eq_pruning_reduced_aim_1} \\
& \mbox{ and } \lim_{\vartheta \to \infty} \limsup_{n \to \infty} \prob^* \Big( \Big| \wtbfN_{n, \varrho}(f) - {\bf N}_{n, \varrho, \vartheta}^\diamond(f)  \Big| > \epsilon/2 \Big) = 0. \label{eq_pruning_reduced_aim_2}
}
Note that \eqref{eq_pruning_reduced_aim_1} and \eqref{eq_pruning_reduced_aim_2} can be proved by modifying the proof of Lemma~3.3 in \cite{bhattacharya:hazra:roy:2016}. Considering the similarity, we postpone the rest of the proof to the subsection~\ref{subsec_pruning_rest_proof}.

\subsection{Proof of Proposition~3.4}

We have noted before that $\prob^*$ is absolutely continuous with respect to $\prob$. To be more precise, we can write down the conditional probability $\prob^*({\sf B}) = q_e^\inv \prob ({\sf B} \cap {\cal S})$ for every ${\sf B} \in {\cal F}$ where $q_e = \exptn(q_e({\bf Y}))$ where $q_e({\bf Y}) = \prob({\cal S} | {\bf Y})$. Define ${\cal S}_{n - \varrho} :=  \{ \mbox{ at least one of $Z_{n - \varrho}$ subtrees is an infinite tree} \}$. It follows immediately that $\mbbo_{\cal S} = \mbbo_{\{ Z_{n - \varrho} > 0\}} \mbbo_{{\cal S}_{n - \varrho}}$. These notations and observations yield the following decomposition of the Laplace functional $\exptn^* (\exp \{ - \int f \dtv \wtbfN_{n, \varrho, \vartheta}\})$
\aln{
q_e^\inv \exptn \Big[ \exp \Big\{ - f \dtv \wtbfN_{n, \varrho, \vartheta} \Big\} \mbbo_{\{Z_{n - \varrho} > 0 \} } \Big] - q_e^\inv \exptn \big[ \exp \big\{  - \int f \dtv \wtbfN_{n, \varrho, \vartheta} \big\} \mbbo_{{\cal S}_{n - \varrho}^c} \big] \label{eq_decomp_first_laplace_functional}
}
where $f \in \ccr$. This decomposition is same as in (4.3) in \cite{bhattacharya:hazra:roy:2016}. The same arguments after the decomposition in the aforementioned reference applies to conclude that the second term in \eqref{eq_decomp_first_laplace_functional} is $o(1)$.  Recall that $\wt{\bbt}_{\varrho, i}$ is the $i$-th regularized tree with root ${\sf r}_i$.  Define
\aln{
\wtbfN^{(i)}_{n, \varrho, \vartheta} := \sum_{\uu \in \wt{\bbt}_{\varrho, i} \setminus \{ {\sf r}_i\}} A^{(\vartheta)} (\uu) \bdelta_{B_n^\inv  X'(\uu)}.  \label{eq_defn_pp_reg_subtree}
}
It is clear that conditioned on the environment ${\bf Y}_{0~:~n-1}$ and $|{\sf D}_{n - \varrho}|$, $(\wtbfN_{n, \varrho, \vartheta}^{(i)} : 1 \le i \le |{\sf D}_{n - \varrho}|)$ are i.i.d. point processes and $\wtbfN_{n, \varrho, \vartheta}$ is the superposition of these i.i.d. point processes. Therefore, the first expectation in \eqref{eq_decomp_first_laplace_functional} equals
\aln{
\exptn \Big[ \mbbo_{\{ |{\sf D}_{n - \varrho}| > 0\}} \Big( \exptn \big[ \exp \{ - \wtbfN_{n, \varrho, \vartheta}^{(1)} \}  \big| ~~ {\bf Y}_{0~:~n- 1}\big] \Big)^{|{\sf D}_{n - \varrho}|} \Big]. \label{eq_first_term_after_decomp}
}
It is clear that the genealogical structure of the point process $\wtbfN_{n, \varrho, \vartheta}^{(1)}$ depends on the segment ${\bf Y}_{n - \varrho~:~n -1}$. But the conditional expectation in \eqref{eq_first_term_after_decomp} can not be restricted to this smaller segment ${\bf Y}_{n - \varrho ~: ~n - 1}$ due to $B_n$ ($B_n$ depends on the full sequence ${\bf Y}_{0~: ~n - 1}$). As the conditional law of the displacements and the genealogical structure interacts, the computations of Laplace functional in \cite{bhattacharya:hazra:roy:2016} and \cite{bhattacharya:hazra:roy:2017} do not apply here as it is.

 To write down the $\wtbfN_{n, \varrho, \vartheta}^{(1)}(f)$ explicitly, we need some notations. If $\uu$ is the $j$-th vertex in the $i$-th generation of $\wt{\bbt}_{\varrho, 1}$, then we shall use the pair $(i, j)$ to denote the vertex $\uu$. Define
\aln{
\wt{\mathbb G}_{\varrho, \vartheta} = [0:~\vartheta^\varrho]^{\vartheta + \vartheta^2 + \ldots + \vartheta^\varrho} \mbox{ and } \wt{\sf R}_0 = \bbr^{\vartheta + \vartheta^2 + \ldots + \vartheta^\varrho} \setminus \{ {\bf 0}_{\vartheta + \vartheta^2 + \ldots + \vartheta^\varrho}\}.
}
We shall  use $\wt{\bf A} = (A^{(\vartheta)}(i, j) : 1 \le j \le \vartheta^i; ~1 \le i \le \varrho)$ and $\wt{\bf X} = (X'(i, j): 1 \le j \le \vartheta^i; ~1 \le i \le \varrho)$.  We can use these notations to write down the the conditional expectation in \eqref{eq_first_term_after_decomp} as follows
\aln{
& \big[ 1 - \pi_n^\inv \sum_{\wt{\bf a} \in \wt{ \mathbb{G}}_{\varrho, \vartheta}} \prob \big(  \wt{\bf A} = \wt{\bf a} | {\bf Y}_{n - \varrho~ :  n - 1}  \big) U_n(\wt{\bf a}) \big] \mbox{ where } \label{eq_term_after_first_decomp_rewritten} \\
& U_n(\wt{\bf a}) :=  \int_{\wt{\sf R}_0} \big[ \pi_n \prob \big( B_n^\inv \wt{\bf X} \in \dtv \wt{\bf x} \big| {\bf Y}_{0 :n- 1}  \big)  \big( 1 - \exp \big\{ - \sum_{i = 1}^\varrho \sum_{j = 1}^{\vartheta^i} a_{i, j} f(x_{i, j}) \big\} \big) \big],\label{added1}
}
and $\wt{\bf a} = (a_{i, j} : 1 \le j \le \vartheta^i;~ 1 \le i \le \varrho)$ and $\wt{\bf x} = (x_{i, j} : 1 \le j \le \vartheta^i ; 1 \le j \le \varrho)$. It is very tempting to think that one must use (stochastic version) of $\mbbM_0$-convergence (for $\pi_n \prob(B_n^\inv \wt{\bf X} \in \cdot | {\bf Y}_{0:~n- 1})$) to obtain the limit of $U_n(\wt{\bf a})$ for every $\wt{\bf a} \in \wt{\mathbb{G}}_{\varrho, \vartheta}$. One needs to be very careful here as the quenched law of $\wt{\bf A}$ also involves $n$ through ${\bf Y}_{n - \varrho ~: ~n - 1}$. To summarize, one must separate out the quenched law of $\wt{\bf X}$ and $\wt{\bf A}$ before letting $n \to \infty$. This can be done as the almost sure limit of $U_n(\wt{\bf a})$ (in $\mbbM_0$-topology) turns out to be  deterministic (does not depend on the environment) and ${\bf Y}_{n - \varrho: n - 1} \eqd {\bf Y}_{0: \varrho - 1}$ and the law of ${\bf Y}_{0: \varrho - 1}$. does not involve $n$. One needs some effort to formalize these ideas and some more notation to write down the limit of  \eqref{eq_term_after_first_decomp_rewritten}.  Recall ${\rm PROJ}_\vartheta : \bbr^\bbn \to \bbr^\vartheta$ from Fact~\ref{fact_finite_dimensional_conv_mzero_conv}. Define $\bsymb{\nu}_*^{(\vartheta)} = \bsymb{\nu}^* \circ {\rm PROJ}_\vartheta$. It follows from \eqref{eq_mzero_conv_crucial} combined with Fact~\ref{fact_finite_dimensional_conv_mzero_conv} (via continuous mapping theorem) that
\aln{
\pi_n \prob \big( B_n^\inv \wt{\bf X} \in \cdot | {\bf Y}_{0~: ~n  - 1} \big) \hlconv \tau(\cdot ) := \sum_{i = 1}^\varrho \sum_{j = 1}^{\vartheta^i} \tau_{i, j}(\cdot) \mbox{ almost surely}
}
in the space $\mathbb{M}(\wt{\sf R}_0)$ where ${\sf F}_i = \{ l \vartheta + 1 : 1 \le l \le \vartheta^{i - 1} - 1\}$ for every $i \ge 1$ and
\aln{
\tau_{i, j} = \underbrace{\bdelta_0 \otimes \bdelta_0  \otimes \cdots \otimes \bdelta_0}_{\vartheta^{i - 1} + j - 1} \otimes \bsymb{\nu}_*^{(\vartheta)} \otimes \underbrace{ \bdelta_0 \otimes \bdelta_0 \otimes \cdots \otimes \bdelta_0}_{ \vartheta^{\varrho} + \vartheta^{\varrho - 1} + \ldots + \vartheta  -  j - \vartheta^{i - 1} - \vartheta + 1} \mbox{ for all } j \in {\sf F}_i \mbox{ and } 1 \le i \le \varrho.
}
The sum and the product structure of $\tau$ can be explained via the dependence structure among the displacements $\wt{\bf X}$. It is well known fact that, for a collection of independent random vectors/processes with (jointly) regularly varying tail, exactly one of the random vectors/processes contribute to the large value of the collection with the same probability (see Subsection~4.5.1 in \cite{lindskog:resnick:roy:2014} and remark after Theorem~3.4 in \cite{resnick:roy:2014}).  In our context, this fact can be realized in the following way. Fix $i \in [1~:~ \varrho]$. Note that if $j \in {\sf F}_i$, then $((i, j), (i, j + 1), \ldots, (i, j + \vartheta - 1))$ have the same parent in the $(i - 1)$-th generation. Thus the displacements $(X'(i, k) : j \le k \le j + \vartheta - 1)$ can be large together ($\bsymb{\nu}_*^{(\vartheta)} $ appears in $\tau_{i, j}$) and none of the other displacements (being independent) contributes to the large value ($\bdelta_0$ appears everywhere else). As $((X'(i, k) : j \le k \le j + \vartheta - 1): j \in {\sf F}_i, ~ 1 \le i \le \varrho)$ are i.i.d. random vectors, the sum over $i$ and $j$  appears in the expression of $\tau$. Recall that ${\bf Y}'$ is an independent copy of ${\bf Y}$ and $W$ is the martingale limit associated to the BPRE.  We are now ready to write down the limit of the expectation in \eqref{eq_first_term_after_decomp}.

\begin{propn} \label{propn_laplace_limn}
If assumptions of Theorem~2.3 hold and we let $n \to \infty$, then the expectation in \eqref{eq_first_term_after_decomp} converges to
\aln{
& q_e \exptn^* \Big[ \exp \Big\{ - W \big[ \prod_{j = 0}^{\varrho - 1} \exptn(\xi | {Y}'_j) \big]^\inv \sum_{\wt{\bf a} \in \wt{\mathbb{G}}_{\varrho, \vartheta}} \prob \big( \wt{\bf A} = \wt{\bf a} \big| {\bf Y}'_{0~: ~\varrho - 1} \big) {\rm U}(\wt{\bf a}) \Big\} \Big]  \label{eq_expsn_Laplace_limn}\\
&  \mbox{ where } {\rm U}(\wt{\bf a})  := \sum_{i = 1} \sum_{j \in {\sf F}_i} \int_{\wt{\sf R}_0} \tau_{i, j}(\dtv \wt{\bf x}) \Big( 1 - \exp \Big\{ - \sum_{i'= 1}^\varrho \sum_{j' \in {\sf F}_i} a_{i', j'} f(x_{i', j'}) \Big\} \Big) \nonumber \\
& = \sum_{i = 1}^\varrho \sum_{j \in {\sf F}_i} \int_{\ol{\bbr}^\vartheta \setminus \{{\bf 0}_\vartheta\}} \bsymb{\nu}_*^{(\vartheta)} \big( \dtv (x_1, x_2, \ldots, x_{\vartheta}) \big) \Big( 1 - \exp \Big\{ - \sum_{k = 1}^\vartheta a_{i, j + k - 1} f(x_k) \Big\} \Big) \mbox{ for all } \wt{\bf a} \in \wt{\mathbb{G}}_{\varrho, \vartheta}.\label{added2}
}
\end{propn}

The proof of this proposition is given in the Subsection~\ref{subsec_proof_laplace_limn}. One can show that there exists a point process ${\bf N}_*^{(\vartheta, \varrho)}$ such that $\exptn^* \big( \exp \big\{ - \int f \dtv {\bf N}_*^{(\vartheta, \varrho)} \big\} \big)$ equals the right hand side of \eqref{eq_expsn_Laplace_limn}. Construction of ${\bf N}_*^{(\varrho, \vartheta)}$ is very similar to the construction of $N_*^{(K, B)}$ in \cite{bhattacharya:hazra:roy:2016} (see equation (4.35) at page 201). So we skip that here. Note that the random variable $W$ appears in the exponent inside expectation \eqref{eq_expsn_Laplace_limn} as the limit of $Z_{n - \varrho}/\pi_{n - \varrho}$ and hence independent of ${\bf Y}'$. In the next step, we shall compute the limit of the \eqref{eq_expsn_Laplace_limn} as $\vartheta \to \infty$. Note that the law of $\wt{\bf A}$ depends on $\vartheta$ due to pruning and the measure $\bsymb{\nu}_*^{(\vartheta)}$ also depends on $\vartheta$. In the next lemma, we shall first write down the sum in the exponent inside expectation in \eqref{eq_expsn_Laplace_limn} in terms of the branching sequence. For that, we need to introduce some notations.

\begin{enumerate}
\item Recall that ${\rm Pow}[1~:~k] $ denotes the power set of $[1~:~k]$ for every $k \ge 1$.
\item Conditioned on ${\bf Y}' = ({\bf y}_i : i \ge 0)$,
 $Z_i^{(\vartheta)}$ denotes the size of the $i$-th generation of a BPRE where each particle in the $(i-1)$-th generation reproduces independently according to the law of $B_{1, i}^{(\vartheta)}$ given through
\aln{
\prob \big( B_{1, i}^{(\vartheta)} = k \big) = \begin{cases} y_{i - 1}(k ) & \mbox{ if } 0 \le  k < \vartheta \\ \sum_{j = \vartheta}^\infty y_{i - 1}(j) & \mbox{ if } k = \vartheta. \end{cases}
}
\item Conditioned on ${\bf Y}'$, $Z_i^{(\vartheta, +)}$ denotes the random variable $Z_i^{(\vartheta)}$ conditioned to stay positive.
\item Conditioned on ${\bf Y}'$, $(Z_i^{(\vartheta, +, k)} : k \ge 1)$ denotes the collection of independent copies of the random variable $Z_i^{(\vartheta, + )}$.
\item Conditioned on ${\bf Y}'_{0~: ~\varrho - 1}$, $\widetilde{Z}_{\varrho - i}^{(\vartheta)}$ denotes the law of the number of descendants in the $\varrho$-th generation of a particle in the $i$-th generation of the BPRE where each particle in the $j$-th generation reproduces independently according to the law of $B_{1,j}^{(\vartheta)}$ for all $j = 1, 2, \ldots, \varrho$.
\end{enumerate}

The following lemma can be derived adapting derivation of (4.24) from (4.21) in \cite{bhattacharya:hazra:roy:2017} and so, the proof of the lemma is omitted as well.

\begin{lemma}
If the assumptions in Theorem~2.3 hold, then we have
\aln{
& \sum_{\wt{\bf a} \in \wt{\mathbb{G}}_{\varrho, \vartheta}} \prob \big( \wt{\bf A} = \wt{\bf a} \big| {\bf Y}'_{0~: ~ \varrho - 1}  \big) {\rm U}(\wt{\bf a})   \nonumber \\
& = \sum_{i = 1}^\varrho \int_{\bbr^\varrho \setminus \{ {\bf 0}_\varrho\}} \bsymb{\nu}_*^{(\vartheta)} \big( \dtv (x_1, x_2, \ldots, x_\vartheta) \big)  \exptn \big( Z_{i - 1}^{(\vartheta)} \big| {\bf Y}'_{0~:~ i - 2}  \big) \prob \Big(Z_1^{(\vartheta)} \ge 1 | Y'_{i -1} \Big) \nonumber \\
& \hspace{.5cm} \exptn \Big[ \sum_{{\sf Q} \in {\rm Pow} ([1 ~: Z_1^{(\vartheta, +)}]) \setminus \{\emptyset\} } \big[ \prob ( Z_{\varrho - i}^{(\vartheta)} \ge 1 ~|~ {\bf Y}'_{i ~: ~\varrho - 1} ) \big]^{|{\sf Q}|}  \big[ \prob( Z_{\varrho - i}^{(\vartheta)} = 0 | {\bf Y}'_{i ~: ~\varrho - 1}) \big]^{Z_1^{(\vartheta, +)} - |{\sf Q}|}  \nonumber \\
& \hspace{1cm} \Big( 1 - \exp \Big\{ - \sum_{k \in {\sf Q}} \wt{Z}_{\varrho - i}^{(\vartheta, +, k)} \Big\} \Big) | {\bf Y}'_{i- 1~: ~\varrho - 1} \Big].  \label{eq_expsn_laplace_limn_branching_version}
}
\end{lemma}

 If we let $\vartheta \to \infty$, then we have $Z_1^{(\vartheta, +)} \asconv Z_1^{(+)}$; $\wt{Z}_{\varrho - i}^{(\vartheta)} \asconv \wt{Z}_{\varrho - i}$ and  $\wt{Z}_{\varrho - i}^{(\vartheta, + )} \asconv \wt{Z}_{\varrho - i}^{(+)}$  conditioned on ${\bf Y}'_{i ~: ~\varrho - 1}$; $\exptn \big( Z_{i - 1}^{(\vartheta)} | {\bf Y}'_{0: i - 2} \big) \asconv \exptn(Z_{i - 1} | {\bf Y}'_{0: ~i - 2})$ for all $1 \le i \le \varrho - 1$. Furthermore, if we let $\vartheta \to \infty$,  Fact~\ref{fact_finite_dimensional_conv_mzero_conv} combined with the dominated convergence  theorem yields that the right hand side of \eqref{eq_expsn_laplace_limn_branching_version} converges to
\aln{
& \sum_{i = 1}^\varrho  \int_{\bbr^\bbn \setminus \{{\bf 0}_\infty\}} \bsymb{\nu}_*(\dtv {\bf x}) \exptn \Big[ Z_{i - 1} | {\bf Y}'_{0~:~i - 2} \Big] \prob( Z_1 > 0 | Y'_{i - 1}) \exptn \Big[ \sum_{{\sf Q} \in {\rm Pow}[1 ~: ~ Z_1^{(+)}]} \big[ \prob \big( Z_{\varrho - i} \ge 1 | {\bf Y}'_{i : ~\varrho - 1} \big) \big]^{|{\sf Q}|} \nonumber \\
& \hspace{1cm} \big[ \prob( Z_{\varrho - i} = 0 | {\bf Y}'_{i ~: ~ \varrho - 1}) \big]^{Z_1^{(+)} - |{\sf Q}|} \Big( 1 - \exp \Big\{ - \sum_{k \in {\sf Q} } \wt{Z}_{\varrho - i}^{(+, k)} f(x_k) \Big\} \Big) \Big| ~~ {\bf Y}'_{i - 1~: ~ \varrho - 1} \Big].  \label{eq_expsn_laplace_functional_lim_vartheta}
}

Our next step would be to let $\varrho \to \infty$. It is natural that the sum in the exponent in \eqref{eq_expsn_Laplace_limn} will become a series but the issue here is that the terms in the sum involves $\varrho$. It is not possible to use some rearrangement of terms such that the terms becomes independent of $\varrho$. We stress this point as the derivation of the Laplace functional in \cite{bhattacharya:hazra:roy:2017} can not help to resolve this fundamental issue.  To circumvent this obstacle, we shall use the exchageability property of the environment which will lead us to a nontrivial distributional identity.  Combining \eqref{eq_expsn_Laplace_limn} and \eqref{eq_expsn_laplace_functional_lim_vartheta}; and ignoring $W$ and $\int \bsymb{\nu}_*^{(\vartheta)}(\dtv {\bf x})$; we obtain the following sum
\aln{
& \sum_{i = 1}^{\varrho}  \Big[ \Big( \prod_{j = 0}^{\varrho - 1} \exptn_{Y'_j}(\xi)\Big)^\inv  \exptn [ Z_{i - 1} | {\bf Y}'_{0~:~i - 2} ]  \Big] \prob( Z_1 > 0 | Y'_{i - 1}) \exptn \Big[ \sum_{{\sf Q} \in {\rm Pow}[1 ~: ~ Z_1^{(+)}]} \big[ \prob \big( Z_{\varrho - i} \ge 1 | {\bf Y}'_{i : ~\varrho - 1} \big) \big]^{|{\sf Q}|} \nonumber \\
& \hspace{1cm} \big[ \prob( Z_{\varrho - i} = 0 | {\bf Y}'_{i ~: ~ \varrho - 1}) \big]^{Z_1^{(+)} - |{\sf Q}|} \Big( 1 - \exp \Big\{ - \sum_{k \in {\sf Q} } \wt{Z}_{\varrho - i}^{(+, k)} f(x_k) \Big\} \Big) \Big| ~~ {\bf Y}'_{i - 1~: ~ \varrho - 1} \Big].
}
It is immediate that the first term in braces equals $[\exptn( Z_{\varrho - i + 1} | {\bf Y}'_{i - 1 ~: ~\varrho - 1})]^\inv$ almost surely conditioned on ${\bf Y}'_{0~: ~\varrho - 1}$. Moreover, being i.i.d. the random variables $(Y'_i : 0 \le i \le \varrho - 1 )$ are exchangeable. This fact lead to the following distributional identity ${\bf Y}'_{0 : \varrho - 1} \eqd {\bf Y}'_{\varrho - 1 : 0}$. Recall that BPRE with environment ${\bf Y}'_{\varrho - 1 ~:~ 0}$ stands for the BPRE where conditioned on ${\bf Y}'_{0: \varrho - 1}$, each particle at the $i$-th generation reproduces independently according to the law $Y'_{\varrho - i + 1}$ for all $i \in [0~:~\varrho - 1]$. Using the fact that  products and the sums are permutation invariant operations and exchangeability of the environment, we have following distributional identity
\aln{
& \sum_{i= 1}^\varrho \big[ \exptn( Z_{\varrho - i + 1} | {\bf Y}'_{i - 1 : \varrho - 1}) \big]^\inv \prob(Z_1 \ge 1 | {\bf Y}'_{i - 1})  \exptn \Big[ \sum_{{\sf Q} \in {\rm Pow}([1 ~: ~ Z_1^{(+)}]) \setminus \{ \emptyset\}}  \big[ \prob( Z_{\varrho - i} \ge 1 | {\bf Y}'_{i : \varrho - 1}) \big]^{|{\sf Q}|} \nonumber \\
& \hspace{.5 cm} \big[ \prob(Z_{\varrho - i} = 0 | {\bf Y}'_{i : \varrho - 1}) \big]^{Z_1^{(+)} - |{\sf Q}|}  \Big( 1 - \exp \big\{ - \sum_{k \in {\sf Q}} \wt{Z}_{\varrho - i}^{(+, k)} f(x_k) \big\} \Big) \Big| {\bf Y}'_{i-1~:~ \varrho - 1}\Big] \nonumber \\
& \eqd \sum_{i = 1}^\varrho \big[ \exptn(Z_{\varrho - i + 1} | {\bf Y}'_{\varrho - i ~:~ 0}) \big]^\inv \prob(Z_1 \ge 1 | Y'_{\varrho -i}) \exptn \Big[ \sum_{{\sf Q} \in {\rm Pow}([1 : Z_1^{(+)}]) \setminus \{ \emptyset \}}  [\prob( Z_{\varrho - i} \ge 1 | {\bf Y}'_{\varrho - i - 1 : 0} )]^{|{\sf Q}|} \nonumber \\
& \hspace{.5cm}  [ \prob( Z_{\varrho - i} = 0 | {\bf Y}'_{\varrho - i - 1: 0})]^{Z_1^{(+)} - |{\sf Q}|} \Big( 1 - \exp \big\{ \sum_{k \in {\sf Q}} \wt{Z}_{\varrho - i}^{(+, k)} f(x_k) \big\} \Big) \Big| {\bf Y}'_{\varrho - i~:~0} \Big]. \label{eq_laplace_identity_exponent_exchangeability}
}
We refer to this identity as the {\it environment seen by the extreme positions}. We can now use the rearrangement of the terms to have following expression for the sum in the right hand side of \eqref{eq_laplace_identity_exponent_exchangeability}
\aln{
& 1 + \sum_{i = 1}^{\varrho - 1} [\exptn(Z_{i + 1} | {\bf Y}'_{i~:~ 0})]^\inv \prob(Z_1 \ge 1 | Y'_i) \exptn \Big[ \sum_{{\sf Q} \in {\rm Pow}([1 : ~Z_1^{(+)}])} [\prob(Z_i \ge 1 | {\bf Y}'_{i-1: 0})]^{|{\sf Q}|}  \nonumber \\
& \hspace{0.5cm} [\prob(Z_i = 0 | {\bf Y}'_{i-1 : 0})]^{Z_1^{(+)} - |{\sf Q}|} \Big( 1 - \exp \big\{ - \sum_{k \in {\sf Q}} \wt{Z}_i^{(+, k)} f(x_k) \big\}\Big) \Big| {\bf Y}'_{i : 0} \Big]. \label{eq_final_exponent_rearrangement}
}
It is easy to check that the right hand side of \eqref{eq_final_exponent_rearrangement} is bounded from above by $1 + \sum_{1 \le i \le \varrho-1}$ $ [\exptn (Z_{i + 1} | {\bf Y}'_{i ~:~ 0}  )]^\inv \le  1 + \sum_{i \ge 1} [ \exptn(Z_{i + 1} | {\bf Y}'_{i :0})]^\inv < \infty$ almost surely due to Remark~1.3. Therefore, we can let $\varrho \to \infty$ in \eqref{eq_expsn_laplace_functional_lim_vartheta}  and use the dominated convergence theorem to obtain the limit which is given in (2.9) in the aforementioned reference.

To conclude that ${\bf N}_n$ converges weakly to ${\bf N}_*$, one needs to derive the Laplace functional of ${\bf N}_*$ and check whether it equals to the limiting Laplace functional or not. This can be achieved by guessing an appropriate {\it marked point process} and writing the Laplace  functional of ${\bf N}_*$ in terms of the Laplace functional of the marked point process. As the derivation is very similar to the derivation of the Laplace functional of $N_*^{(K, B)}$ in \cite{bhattacharya:hazra:roy:2017} (see equation (4.35) at page~201), we decided to provide only the brief derivation in the Subsection~\ref{subsdubsec_appendix_characterization_N_star}.

\section{Further proofs} \label{sec_further_proofs}



\subsubsection{Rest of the proof of Proposition~3.3 } \label{subsec_pruning_rest_proof}
{\bf Proof of \eqref{eq_pruning_reduced_aim_2}.} Note that
\aln{
\big| \wtbfN_{n, \varrho}(f) - \wtbfN_{n, \varrho, \vartheta}^{\diamond}(f) \big| \le \norm{f} \sum_{i =1}^{|{\sf D}_{n - \varrho}|} \sum_{\sfu \in \wh{\sf T}_{\varrho, i}} \big[ A(\sfu) - A^{(\vartheta)}(\sfu) \big] \mbbo(|\Xu| > B_n \delta)
}
almost surely. Recalling ${\sf J}_{n, \delta}(\epsilon)$ introduced in \eqref{eq_old_defn_Jn_eta_wt_vep}, we can derive the following upper bound for the probability in \eqref{eq_pruning_reduced_aim_1}
\aln{
& \exptn^* \Big[ \mbbo \big( {\sf J}_{n, \delta}(\epsilon) \big) \prob^* \Big( \norm{f} \sum_{ i =1}^{|{\sf D}_{n - \varrho}|} \sum_{\sfu \in \wh{\sf T}_{\varrho, i}} \big[ A(\sfu) - A^\diamond(\sfu) \big] \mbbo(|X| > B_n \delta) > \varepsilon | {\bf Y}_{0:~n - 1} \Big) \Big] +  \prob^* \big( {\sf J}_{n, \delta}^c(\epsilon) \big) \nonumber \\
& \le \varepsilon^\inv \norm{f} \exptn^* \Big[ \mbbo({\sf J}_{n, \delta}(\epsilon)) \exptn^* \Big( \sum_{i =1}^{|{\sf D}_{n - \varrho}|} \sum_{\sfu \in \wh{\sf T}_{\varrho, i}} \big[ A(\sfu) - A^\diamond(\sfu) \big]  \mbbo(|\Xu| > B_n \delta) \big| {\bf Y}_{0:~n-1} \Big) \Big] + o(1) \label{eq_comp_cut_mark_markov_bound}
}
using Markov inequality.  We now treat the conditional expectation in  \eqref{eq_comp_cut_mark_markov_bound}. The first observation is that ${\sf D}_{n - \varrho}$ and the displacements indexed by the forest $(\wh{\sf T}_{\varrho, i} : 1 \le i \le |{\sf D}_{n - \varrho}|)$ are independent conditioned on ${\bf Y}_{0:~n-1}$. This observation leads to the following expression for the conditional expectation
\aln{
\exptn^* \Big( |{\sf D}_{n - \varrho}| \big| {\bf Y}_{0:~n - 1} \Big) \sum_{i =1}^\varrho \exptn^* \Big( \sum_{\sfu \in \wh{\sf T}_{\varrho, 1}} \big[ A(\sfu) - A^\diamond(\sfu) \big] \mbbo(|\Xu| > B_n \delta) \Big| {\bf Y}_{0: ~n- 1} \Big) \label{eq_comp_cut_marked_cond_exptn_equality}
}
almost surely. We now note that the branching mechanism and $(\mbbo(|\Xu| > B_n \delta) : \sfu \in \wh{\sf T}_{\varrho, 1})$ are independently distributed conditioned on the environment ${\bf Y}_{0:~n - 1}$. Using the fact that $(\mbbo(|\Xu| > B_n \delta) : \sfu \in \wh{\sf T}_{\varrho, 1})$ are identically distributed we can obtain the following upper bound
\aln{
\prob(|X_1| > B_n \delta | {\bf Y}_{0:~n-1}) \sum_{i =1}^\varrho \exptn^* \Big( \sum_{\sfu \in \wh{\sf D}_i^{(1)}} \big[ A(\sfu) - A^\diamond(\sfu) \big] \big| {\bf Y}_{0:~n-1} \Big) \label{eq_comp_cut_marked_ident_disp}
}
almost surely for the conditional expectation inside the sum in \eqref{eq_comp_cut_marked_cond_exptn_equality}. We now note that conditioned on ${\bf Y}_{0:~n - 1}$, $\wh{\sf D}_i \eqd Z_i$ are independent of $A(\sfu)$ and $A^\diamond(\sfu)$ where $A(\sfu) \eqd {\wt Z}_{\varrho - i}$ and $A^\diamond(\sfu) = \wt{ Z}^{(\vartheta)}_{\varrho - i}$ if $\sfu \in \wh{\sf D}_i$, where
$\wt{ Z}^{(\vartheta)}_i$ is the number of descendants in the $n$-th generation of a particle
in the $(n - i)$-th generation conditioned on the environment ${\bf Y}_{0:~n- 1}$ with progeny random variable $Z_1^{(\vartheta)} = Z_1 \wedge \vartheta$.
These observations together yield following expression
\aln{
\sum_{i =1}^\varrho  \exptn^* \Big( \wt{ Z}_{\varrho - i}  - \wt{ Z}^{(\vartheta)}_{\varrho - i} \big| {\bf Y}_{0:~n- 1}\Big)  \label{eq_comp_cut_marked_au}
}
for the conditional expectation in \eqref{eq_comp_cut_marked_ident_disp}. Combining the expressions obtained in \eqref{eq_comp_cut_marked_cond_exptn_equality}- 
\eqref{eq_comp_cut_marked_au} and using definition \eqref{eq_old_defn_Jn_eta_wt_vep} of the set ${\sf J}_{n, \delta}(\epsilon)$, we have the the following upper bound for the conditional expectation in \eqref{eq_comp_cut_mark_markov_bound}
\aln{
& \mbbo({\sf J}_{n, \delta}(\epsilon))\Big[ \pi_n \prob( |X_1| > B_n \delta | {\bf Y}_{0:~n - 1}) \Big] \frac{\pi_{n - \varrho}}{\pi_n} \sum_{i =1}^\varrho \Big( \prod_{j = 0}^{i -1}
\exptn(\xi|Y_{n - \varrho + j}) \Big) \nonumber \\
& \le  (\delta^{- \alpha} + \epsilon) \sum_{i =1}^\varrho \Big( \prod_{i = n - \varrho}^{n - 1} \exptn(\xi|Y_i) \Big)^\inv \exptn^* \Big( \wt{ Z}_{\varrho - i} - \wt{ Z}_{\varrho - i}^{(\vartheta)} \big| {\bf Y}_{0:~n - 1} \Big)
 \nonumber \\
& = (\delta^{- \alpha} + \epsilon) \sum_{i =1}^\varrho  \Big( \prod_{j = i}^{\varrho - 1} \exptn(\xi|Y_{n - \varrho + j}) \Big)^\inv \Big[ \Big( \prod_{j'= i}^{\varrho - 1} \exptn(\xi|Y_{n - \varrho + j'})  \Big) - \Big( \prod_{j'= i}^{\varrho - 1} \exptn(\xi \wedge \vartheta|Y_{n - \varrho + j'}) \Big)  \Big] \nonumber \\
& = (\delta^{- \alpha} + \epsilon) \sum_{ i =1}^{\varrho} \Big[ 1 -  \prod_{j =i}^{\varrho - 1} \frac{ \exptn(\xi \wedge \vartheta|Y_{n - \varrho - j})}{\exptn(\xi|Y_{n - \varrho + j})}  \Big]
}
almost surely. Plugging this upper bound back in \eqref{eq_comp_cut_mark_markov_bound}, we obtain
\aln{
\varepsilon^\inv \norm{f}(\delta^{- \alpha} + \epsilon) \sum_{i =1}^\varrho \Big( 1 - \exptn \Big[ \prod_{j =0}^{\varrho - i - 1} \frac{\exptn(\xi\wedge \vartheta|Y_j)}{\exptn(\xi|Y_j)} \Big] \Big) \label{eq_comp_cut_marked_final_upper_bound}
}
using $(Y_{n - \varrho + i + j} : 0 \le j \le \varrho - i -1) \eqd (Y_j : 0 \le j \le \varrho - i)$ for every $ 1 \le i \le \varrho$ due to stationarity of the environment. The upper bound then is independent of $n$. Note that the sum and the product in \eqref{eq_comp_cut_marked_final_upper_bound} is finite for fixed $\varrho \ge 1$. Also note that the product is bounded by one  almost surely and so we can apply the dominated convergence theorem as $\vartheta \nearrow \infty$. Hence the upper bound in \eqref{eq_comp_cut_marked_final_upper_bound} converges to zero as $\vartheta \nearrow \infty$.\\

\noindent{\bf Proof of \eqref{eq_pruning_reduced_aim_1}.} We start with the observation ${\sf D}_j^{(\diamond, i)} \subset \wh{\sf D}_j^{(i)}$ and $A^{\diamond}(\sfu) = A^{(\vartheta)}(\sfu)$ for all $\sfu \in {\sf D}^{(\diamond, i)}_j$ for all $i \ge 1$ and $j \ge 1$. We then observe that
\aln{
& {\bf N}^{\diamond}_{n, \varrho, \vartheta}(f) - \wtbfN_{n, \varrho, \vartheta}(f)  \nonumber \\
&\quad = \sum_{i =1}^{|{\sf D}_{n - \varrho}|} \sum_{j =1}^\varrho \sum_{\sfu \in \wh{\sf D}_j^{(i)} \setminus {\sf D}^{(\diamond, i)}_j} A^{\diamond}(\sfu) f( B_n^\inv \Xu) \nonumber \\
& \quad\le \norm{f} \sum_{i =1}^{|{\sf D}_{n - \varrho}|} \sum_{j =1}^\varrho \sum_{\sfu \in \wh{\sf D}^{(i)}_j \setminus {\sf D}^{(\diamond, i)}_j} A^{\diamond}(\sfu)  \mbbo(|\Xu| > B_n \delta)
}
almost surely. Then we can follow the similar steps like in the proof of \eqref{eq_pruning_reduced_aim_2} and derive the following upper bound for the probability in \eqref{eq_pruning_reduced_aim_1}
\aln{
& \norm{f} \varepsilon^\inv \exptn^* \Big[ \mbbo({\sf J}_{n, \delta}(\epsilon)) \exptn^* \Big( \sum_{i =1}^{|{\sf D}_{n - \varrho}|} \sum_{j =1}^\varrho \sum_{\sfu \in \wh{\sf D}_j^{(i)} \setminus {\sf D}_j^{(\diamond, i)}} A^\diamond(\sfu) \nonumber\\
& \hspace{2cm}  \mbbo( |\Xu| > B_n \delta) \big| {\bf Y}_{0:~n - 1} \Big) \Big] + \prob^*({\sf J}_{n, \delta}^c(\epsilon)),
}
where ${\sf J}_{n, \delta}(\epsilon)$ is defined in \eqref{eq_old_defn_Jn_eta_wt_vep}
and $\prob({\sf J}_{n, \delta}^c(\epsilon)) = o(1)$. We shall use $\exptn(\xi | Y_i)$ to denote the mean of the (random) distribution $Y_i$ for $i \ge 0$.
Using similar arguments like in the proof of \eqref{eq_aim_old_event}
and stationarity of the environment we can derive the following upper bound
\aln{
& \varepsilon^\inv \norm{f}(\delta^{- \alpha} + \epsilon) \exptn^* \Big[ \Big( \prod_{j = n - \varrho}^{n - 1} \exptn(\xi|Y_j) \Big)^\inv \sum_{j =1}^\varrho \exptn \Big( \sum_{\sfu \in \wh{\sf D}_j^{(i)} \setminus {\sf D}_j^{(\diamond, i)}} A^\diamond(\sfu) \big| {\bf Y}_{0:~n - 1} \Big) \Big] \nonumber \\
& =\varepsilon^\inv \norm{f} (\delta^{- \alpha} + \epsilon) \exptn^* \Big[ \Big( \prod_{j = n - \varrho}^{n - 1} \exptn(\xi|Y_j) \Big)^\inv \sum_{j =1}^\varrho \exptn \big( Z_j - Z_j^{(\vartheta)} | {\bf Y}_{n - \varrho:~n- 1} \big) \exptn(\wt{Z}^{(\vartheta)}_{\varrho - j} | {\bf Y}_{n - \varrho ~:~n - 1}) \Big] \nonumber \\
& \le \varepsilon^\inv \norm{f} (\delta^{- \alpha} + \epsilon) \exptn^* \Big[ \Big( \prod_{j = n- \varrho}^{n - 1} \exptn(\xi|Y_j) \Big)^\inv \sum_{j =1}^\varrho \exptn \big( Z_j - Z_j^{(\vartheta)} \big| {\bf Y}_{0:~n- 1} \big) \Big] \nonumber \\
& = \varepsilon^\inv \norm{f} (\delta^{- \alpha} + \epsilon) \sum_{i =1}^\varrho \Big[ 1 - \exptn \Big( \prod_{j =0}^{i -1} \frac{\exptn(\xi \wedge \vartheta|Y_{n - \varrho + j})}{ \exptn(\xi|Y_{n - \varrho + j})} \Big) \Big] \nonumber \\
& = \varepsilon^\inv \norm{f} (\delta^{- \alpha} + \epsilon) \sum_{i =1}^\varrho \Big[ 1 - \exptn \Big( \prod_{j =0}^{i -1} \frac{\exptn(\xi \wedge \vartheta|Y_j)}{\exptn(\xi|Y_j)} \Big) \Big]. \label{eq_comp_marked_prun_final_bound}
}
The upper bound obtained in \eqref{eq_comp_marked_prun_final_bound} is independent of $n$. The sum is finite as $\varrho$ is finite. We can let $\vartheta \nearrow \infty$ and apply dominated convergence theorem to conclude that each term inside the sum converges to zero. Hence the proof is complete.

\subsection{Proof of Proposition~\ref{propn_laplace_limn}} \label{subsec_proof_laplace_limn}

Fix $\epsilon > 0$ and define
\aln{
{\sf H}_{n, \epsilon} := \big\{ \max_{\wt{\bf a} \in \wt{S}_{\varrho, \vartheta}} \big| U_n(\wt{\bf a}) - {\rm U}(\wt{\bf a}) \big| < \epsilon \big\}.
}
From \eqref{added1} and \eqref{added2} it follows that
$U_n(\wt{\bf a}) \asconv {\rm U}(\wt{\bf a})$ as $n \rightarrow+\infty$. Further, $|\wt{\mathbb S}_{\varrho, \vartheta}| < \infty$
and thus
%
\aln{
\lim_{n \to \infty} \prob \big( {\sf H}_{n, \epsilon}^c \big) = 0.
\label{ccweak}}
Therefore, it is enough to compute the limit of Laplace functional on the event ${\sf H}_{n, \epsilon}$
\aln{
\lim_{n \to \infty} \exptn \Big[ \mbbo({\sf H}_{n, \epsilon}) \mbbo(|{\sf D}_{n - \varrho}| > 0) \Big( \exptn \big[ \exp \big\{ - {\bf N}_{n, \varrho, \vartheta}^{(1)}(f) \big\} \big]  \Big| {\bf Y}_{0~:~n - 1}\Big)^{|{\sf D}_{n - \varrho}|} \Big]  \label{eq_aim_limn_laplace}
}
when $\epsilon$ is arbitrarily small.  Observe that 
\aln{
& \mbbo({\sf H}_{n, \epsilon}) \Big[ 1 - \pi_n^\inv \sum_{\wt{\bf a} \in \wt{\mathbb S}_{\varrho, \vartheta}} \prob \big( \wt{\bf A} = \wt{\bf a} | {\bf Y}_{0~:~n - 1} \big) U_n(\wt{\bf a}) \Big]^{|{\sf D}_{n - \varrho}|} \mbbo(|{\sf D}_{n - \varrho}| > 0) \nonumber \\
& \le
\Big[ 1 - \pi_{n - \varrho} \Big( \prod_{j = n - \varrho}^{n - 1} \exptn_{Y_j}(\xi) \Big)^\inv \sum_{\wt{\bf a} \in \wt{\mathbb S}_{\varrho, \vartheta}} \prob \big( \wt{\bf A} = \wt{\bf a} | {\bf Y}'_{n- \varrho ~: ~n - 1} \big)  \nonumber \\
& \hspace{2cm} \Big({\rm U}(\wt{\bf a}) - \epsilon \Big) \Big]^{\pi_{n - \varrho} (|{\sf D}_{n - \varrho}|/ \pi_{n - \varrho})} \label{eq_propn_cutting_decomp_reg_dep}
}
almost surely. We first note that $\pi_{n - \varrho}$ and ${\sf D}_{n - \varrho}$ are independent of ${\bf Y}_{n - \varrho~:~n - 1}$ as the environment random variables are independently distributed. To use the i.i.d. structure of the environment, consider an independent copy ${\bf Y}'= (Y'_i : i \ge 0)$ of ${\bf Y}$ and then we can see that the expression in \eqref{eq_propn_cutting_decomp_reg_dep} equals in distribution to the following expression
\aln{
& \Big[ 1 - \pi_{n - \varrho}^\inv \Big( \prod_{j =0}^{\varrho - 1} \exptn_{Y'_j}(\xi) \Big)^\inv \sum_{\wt{\bf a} \in \wt{\mathbb S}_{\varrho, \vartheta}} \prob \big( \wt{\bf A} = \wt{\bf a} \big| {\bf Y}'_{0:~\varrho - 1} \big) \nonumber \\
& \hspace{1cm} \Big( {\rm U}(\wt{\bf a}) - \epsilon  \Big) \Big]^{\pi_{n - \varrho} (|{\sf D}_{n - \varrho}|/ \pi_{n - \varrho})} \mbbo(|{\sf D}_{n - \varrho}| > 0). \label{eq_ind_cutting_decomp_reg_dep}
}
We know that $\mbbo(|{\sf D}_{n - \varrho}| > 0) \asconv \mbbo({\cal S})$ as $n \to \infty$ and
$\pi_{n - \varrho} \asconv \infty$. Therefore, on the event ${\cal S}$, the expression obtained in \eqref{eq_ind_cutting_decomp_reg_dep} converges almost surely to
\aln{
\exp \Big\{ - W \Big( \prod_{j =0}^{\varrho - 1} \exptn_{Y'_j}(\xi) \Big)^\inv \sum_{\wt{\bf a} \in \wt{\mathbb S}_{\varrho, \vartheta}} \prob \big( \wt{\bf A} = \wt{\bf a} | {\bf Y}'_{0:~\varrho - 1} \big)  \Big( {\rm U}(\wt{\bf a}) - \epsilon \Big) \Big\}. \label{eq_aslim_n_propn_reg_dep}
}
Then the dominated convergence theorem implies that
\aln{
& \limsup_{n \to \infty} \exptn \Big[ \mbbo({\sf H}_{n, \epsilon}) \mbbo(|{\sf D}_{n - \varrho}| > 0) \Big( \exptn \big[e^{ - \wt{\bf N}_{n, \varrho, \vartheta}^{(1)}(f)} | {\bf Y}_{0:~n- 1}\big] \Big)^{|{\sf D}_{n - \varrho}|} \Big] \nonumber \\
& \le \exptn \Big[ \mbbo({\cal S}) \exp \Big\{ - W \Big( \prod_{j =0}^{\varrho - 1} \exptn(\xi|Y'_j) \Big)^\inv   \sum_{\wt{\bf a} \in \wt{\mathbb S}_{\varrho, \vartheta}} \prob \big( \wt{\bf A} = \wt{\bf a} | {\bf Y}'_{0:~\varrho - 1} \big)  \Big( {\rm U}(\wt{\bf a}) - \epsilon \Big)\Big\} \Big]. \label{eq_limsup_n_propn_reg_dep}
}
We can follow a similar path starting from the derivation in \eqref{eq_propn_cutting_decomp_reg_dep} that
\aln{
& \liminf_{n \to \infty} \exptn \Big[ \mbbo({\sf H}_{n, \epsilon}) \mbbo(|{\sf D}_{n - \varrho}| > 0) \Big( \exptn \big[e^{ - \wt{\bf N}_{n, \varrho, \vartheta}^{(1)}(f)} | {\bf Y}_{0:~n- 1}\big] \Big)^{|{\sf D}_{n - \varrho}|} \Big] \nonumber \\
& \ge \exptn \Big[ \mbbo({\cal S}) \exp \Big\{ - W \Big( \prod_{j =0}^{\varrho - 1} \exptn(\xi|Y'_j) \Big)^\inv   \sum_{\wt{\bf a} \in \wt{\cal S}_{\varrho, \vartheta}} \prob \big( \wt{\bf A} = \wt{\bf a} | {\bf Y}'_{0:~\varrho - 1} \big)  \Big( {\rm U}(\wt{\bf a}) + \epsilon \Big)\Big\} \Big]. \label{eq_liminf_n_propn_reg_dep}
}
The proof of this lemma follows from \eqref{eq_limsup_n_propn_reg_dep} and \eqref{eq_liminf_n_propn_reg_dep} by letting $\epsilon \to 0$.

\subsection{Characterization of the weak limit ${\bf N}_*$} \label{subsdubsec_appendix_characterization_N_star}

Define the marked Cox process
\aln{
{\bf N}_*^\marked : = \sum_{l = 1}^\infty \bdelta_{(V_l, {\bf R}^{(l)}, [{\rm C}_1({\bf Y}') W]^{1/\alpha} \bsymb{\zeta}^{(l)})}
}
on the state space $\mathscr{X} = \cup_{v = 1}^\infty \{ \{ v\} \times \bbn_{\#}^v \times \bbr_{{\bf 0}}^\bbn\}$.  To understand the behavior of the point ${\mathscr X}$-valued point process, we need to look at the class of all bounded continuous functions $h : \mathscr{X} \to \bbr$ with compact support. Note that $\mathscr{X}$ can be embedded in the larger metric space ${\mathscr X}_* = \bbn \times {\{0\} \cup\bbn}^{\bbn} \times \bbr^{\bbn}_{\bf 0}$ adding appropriate number of $0$'s. The metric ${\bf d}_*$ on $\mathscr{X}_*$ is given by ${\bf d}_*({\bf s}, {\bf t}) = |s_1 - t_1| + \sum_{i = 1}^\infty 2^{- i} (|s_{2, i} - t_{2,i}| \wedge 1) + \sum_{i = 1}^\infty 2^{-i} (|s_{3,i} - t_{3,i}| \wedge 1)$ where ${\bf s} = (s_1, (s_{2, i} : i \ge 1), (s_{3,i} : i \ge 1))$ and ${\bf t} = (t_1, (t_{2, i} : i \ge 1), (t_{3,i} : i \ge 1))$.   Conditioned on ${\bf Y}'$ and $W$, ${\bf N}_*^\marked$ is a marked Poisson process with mean measure $\exptn [{\bf N}_*^\marked ( \{v\} \times {\sf K} \times {\sf H})]$ given by
\aln{
\bsymb{\nu}_*\big( [{\rm C}_1({\bf Y}') W]^{- 1/\alpha} . {\sf H} \big) \otimes \prob \big[ (V, {\bf R}) \in \{ v \} \times {\sf K} | {\bf Y}'  \big]
}
where $v \in \bbn$, ${\sf K} \subset \bbn_{\#}^v$ and ${\sf H} \in {\cal B}(\bbr_{{\bf 0}}^\bbn)$ such that ${\bf 0}_\infty \notin {\rm cl}({\sf H})$. Consider $h \in C_{c}^+({\mathscr X})$. We can use Proposition~3.8 in \cite{resnick:1987} and homogeneity property of $\bsymb{\nu}_*$ (see (2.1)) to obtain the following expression for the Laplace functional  $\exptn^*( \exp \{ - {\bf N}_*^\marked(h)\}) = \exptn^* [\exp \{ - \sum_{l = 1}^\infty h(V_l, {\bf R}^{(l)}, \bsymb{\zeta}_l)\}]$
\aln{
\exptn^* \Big( \exp \Big\{ - {\rm C}_1({\bf Y}') W \int_{\bbr_{\bf 0}^\bbn} \bsymb{\nu}_*(\dtv x) \exptn \Big( 1 - \exp \big\{ - h(V, {\bf R}, {\bf x}) \big\} \Big|~ {\bf Y}'  \Big) \Big\}\Big). \label{eq_laplace_functional_marked_pp}
}
Consider $f \in \ccr$. Consider $h_0: \mathscr{X} \to \bbr$ such that for all ${\bf x} = (v, (x_{1, i} : 1 \le i \le v), (x_{2, i} : i \ge 1)) \in \mathscr{X}$, $h_0({\bf x}) = \sum_{j = 1}^v x_{1, j} f(x_{2,j})$. It is easy to check that $h_0 \in C_c^+(\mathscr{X})$ and ${\bf N}_*(f) = {\bf N}_*^{\marked}(h_0)$ almost surely.  It follows immediately from \eqref{eq_laplace_functional_marked_pp} that the conditional expectation in the exponent of  $\exptn^*(\exp\{ - {\bf N}_*(f) \})$ equals $\exptn ( 1 - \exp \{ - \sum_{i = 1}^V R_i f(x_i) \} |{\bf Y}' )$.
We can use the joint conditional p.m.f of $(V, {\bf R})$ in (2.6) to conclude that $\exptn( \exp \{ - {\bf N}_*(f)\})$ equals the expression in (2.9) in the aforementioned reference. Hence, the series representation of the weak limit ${\bf N}_*$ follows.
\end{document}